# NEW ENHANCED CHAOTIC NUMBER GENERATORS

**René Lozi**

**Laboratoire J.A. Dieudonné – UMR du CNRS N° 6621,
University of Nice-Sophia-Antipolis, Parc Valrose, 06108 Nice Cedex 02, France**

**and**

**Institut Universitaire de Formation des Maîtres Célestin Freinet-académie de Nice,
89 avenue George V, 06046 Nice Cedex 1, France,  rlozi@unice.fr**

## Abstract

We introduce new families of enhanced chaotic number generators in order to compute very fast long series of pseudorandom numbers.

Generation of random or pseudorandom numbers in nowadays is a key feature for industrial mathematics. Pseudorandom or chaotic numbers are used in many areas of contemporary technology such as modern communication systems and engineering applications. Everything we do to achieve privacy and security in the computer age depends on random numbers. More and more patents using discrete mappings for this purpose are obtained by researchers of discrete dynamical systems. The idea of construction of chaotic pseudorandom number generators (CPRNG) applying discrete chaotic dynamical systems intrinsically exploits the property of extreme sensitivity of trajectories to small changes of initial conditions, since the generated bits are associated with trajectories in an appropriate way.

Recently we have brought to light that the very very weak coupling of many symmetric tent maps (and logistic maps) can generate chaotic numbers in good conditions. Moreover these numbers are equally distributed over a given finite interval. Numerical computations show that this distribution is obtained with a very good approximation. They have also the property that the length of the periods of the numerically observed orbits is very large. However chaotic numbers are not pseudo-random numbers because the plot of the couples of iterated points $(x_n, x_{n+1})$ in the phase plane shows up the symmetric tent map used to generate them.

We introduce now new families of enhanced chaotic number generators in order to hide the generating function. The key feature of these enhanced chaotic number generators is that they use chaotic numbers themselves in order to sample chaotic subsequences of chaotic numbers.

We explore numerically the properties of these new families and underline their very high qualities and usefulness as CPRNG when series are computed up to $10^{13}$ iterations.

# 1. Introduction

Generation of random or pseudorandom numbers nowadays is a key feature for industrial mathematics. Pseudorandom or chaotic numbers are used in many areas of contemporary technology such as modern communication systems and engineering applications.

As an example, Monte Carlo methods are a widely used class of computational algorithms for simulating the behavior of various physical and mathematical systems, and for other computations [1]. They are distinguished from other simulation methods by being nondeterministic in some manner, usually by using random numbers (or, more often, pseudo-random numbers). In mathematics, a Monte Carlo algorithm is a numerical Monte Carlo method used to find solutions to mathematical problems (which may have many variables) that cannot easily be solved, for example, by integral calculus, or other numerical methods. For many types of problems, its efficiency relative to other numerical methods increases as the dimension of the problem increases.

Indeed, deterministic methods of numerical integration operate by taking a number of evenly spaced samples from a function. In general, this works very well for functions of one variable. However, for functions of vectors, deterministic quadrature methods can be very inefficient. To numerically integrate a function of a two-dimensional vector, equally spaced grid points over a two-dimensional surface are required. For instance a 10x10 grid requires 100 points. If the vector has 100 dimensions, the same spacing on the grid would require $10^{100}$ (one **googol**) points—that's far too many to be computed. In fact 100 dimensions are by no means unreasonable, since in many physical problems, a "dimension" is equivalent to a degree of freedom.

Monte Carlo methods provide a way out of this exponential time-increase. As long as the function in question is reasonably well-behaved, it can be estimated by randomly selecting points in 100-dimensional space, and taking some kind of average of the function values at these points. By the law of large numbers, this method will display $\frac{1}{\sqrt{N}}$ convergence—i.e. quadrupling the number of sampled points will halve the error, regardless of the number of dimensions.

The new families of enhanced chaotic number generators we introduce in this article are able to generate pseudo-random numbers not only in one dimension but also in any dimensions with uncorrelated coordinates as needed in this case.

Monte Carlo methods are also very important in computational physics, physical chemistry, and related applied fields, and have diverse applications from complicated quantum chromodynamics calculations to designing heat shields and aerodynamic forms.

Monte Carlo methods have likewise proven efficient in solving coupled integral differential equations of radiation fields and energy transport, and thus these methods have been used in global illumination computations which produce photorealistic images of virtual 3D models, with applications in video games, go and chess programming, architecture, design, computer generated films, special effects in cinema, business, economics and other fields.

Monte Carlo methods are in addition useful in many areas of computational mathematics, where a lucky choice can find the correct result. A classic example is Rabin's algorithm for primality testing.

Even though widely used Monte Carlo methods are not the only methods which require random or pseudorandom numbers. Random, pseudo-random and chaotic numbers are found very useful for industrial mathematics for a variety of purposes such as generation of cryptographic keys, computer games and some classes of scientific experiments. Everything we do to achieve privacy and security in the computer age depends on random numbers [2, 3].

Computer engineers chose to introduce randomness into computers in the form of *pseudo-random number generators* (PRNG). As the name suggests, pseudo-random numbers are not truly random. Rather, they are computed from a mathematical formula or simply taken from a precalculed list. A lot of research has gone into pseudo-random number theory and modern algorithms for generating them are so good that the numbers look like they were really random.

There are several requirements for a good PRNG and its implementation in a subroutine library. Among them are statistical robustness (uniform distribution of values at the output with no apparent correlations), unpredictability, long period, efficiency, theoretical support (precise prediction of the important properties), portability, and others [4-7]. A number of PRNGs introduced in the last five decades fulfil most of the requirements and are successfully used in simulations (*e.g.* the Mersenne Twister PNRG [8]). Nevertheless, each of them has some weak properties which may - or may not - influence the results.

It has been observed that PRNGs are examples of deterministic chaotic dynamical [9], this explains why a great care must be taken in the implementation of these generators so that no floating-point truncation/rounding errors occurs.

The idea of construction of chaotic pseudorandom number generators (CPRNG) applying discrete chaotic dynamical systems intrinsically exploits the property of extreme sensitivity of trajectories to small changes of initial conditions, since the generated bits are associated with trajectories in an appropriate way [10].

Recently some authors [11] proposed the use of the Arnol'd cat maps as a PNRG (Arnol'd cat map are hyperbolic automorphisms which are represented by 2x2 matrices with integer entries, a unit determinant, and real eigenvalues). Then the properties making a deterministic algorithm suitable to generate a pseudo-random sequence of numbers have been discussed: high value of Kolmogorov-Sinai entropy, high-dimensionality for the parent dynamical system, and very large period of the generated sequence. A multi-dimensional version of the Arnol'd cat maps has also been proposed as a PNRG [12].

This field of research is nowadays in constant development [13-15]. From an industrial point of view, more and more European or United States patents using discrete mapping to produce pseudo-random numbers are obtained by researchers of discrete dynamical systems [16-20].

Recently we have brought to light that the very weak coupling of the symmetric tent map can generate deftly chaotic numbers in one or several dimensions [21]. Moreover these numbers are equally distributed over a given finite interval. Numerical computations show that this distribution is obtained with a very good approximation. They have also the property that the length of the periods of the numerically observed orbits is very large. However chaotic numbers are not pseudo-random numbers because the plot of the couples of iterates ($x_n$, $x_{n+1}$) in the phase plane shows up the symmetric tent map used to generate them.

We introduce now new families of enhanced chaotic number generators in order to hide the generating function. The key feature of these enhanced chaotic number generators is that they use chaotic numbers themselves in order to sample chaotic subsequences of chaotic numbers.

We explore numerically the properties of these new families and underline their very high qualities and usefulness as CPRNG when series are computed up to $10^{13}$ iterations.

The process of chaotic sampling of chaotic sequences which is pivotal for these new families works perfectly in numerical simulation when floating point (or double precision) numbers are handled by computer.

It is noteworthy that the new models of very weakly coupled maps we introduce are more powerful than the usual chaotic maps used to generate chaotic sequences mainly because only additions and multiplications are used in the computation process, no division being required. Moreover the computations are done using floating point or double precision numbers, allowing the use of the powerful Floating Point Unit (FPU) of the modern microprocessors (built by both Intel or AMD companies). In addition a large part of the computations can be parallelized taking advantage of the multicore microprocessors.

In section 2, we recall the basic equations of the coupled symmetric tent maps introduced in [21] and we give more precise numerical results on them.

In Section 3 we introduce new families of enhanced chaotic number generators in order to hide the generating function using chaotic numbers themselves in order to sample chaotic subsequences of chaotic numbers.

We explore numerically the properties of these new families and underline their very high qualities and usefulness as CPRNG when series are computed up to $10^{13}$ iterations.

## 2. Weakly coupled symmetric tent maps as chaotic generators

When a dynamical system is realized on a computer using floating point or double precision numbers, the computation is of a discretization, where finite machine arithmetic replaces continuum state space. For chaotic dynamical systems, the discretization often have collapsing effects to a fixed point or to short cycles.

As an example of such collapsing effects, O. E. Lanford III [22] presents the results of a sampling study in double precision of a discretization of the logistic map which has excellent ergodic properties. However the very long-term behaviour of numerical orbits is, for a substantial fraction of initial points, in flagrant disagreement with the true behaviour of typical orbits of the original smooth logistic map. In [23] the authors study an example of random map where the components have absolutely continuous invariant measures (*acim*), but where computer experiments reveal the surprising fact that all orbits eventually fall into a stable periodic orbit. They explain why the *acim* cannot be observed experimentally.

In order to preserve the chaotic properties of the continuous models in numerical experiments we have introduced the very weak coupling of one-dimensional dynamical systems which are noteworthy.

### 2.1 multi-dimensional coupling

#### 2.1.1 Two-dimensional coupled map

We recall first the basic equation of the coupled symmetric tent maps introduced in [21]. In this paper we will consider only the symmetric tent map defined by

$$f_a(x) = 1 - a|x| \qquad (2.1)$$

with the value $a = 2$, later denoted simply as $f$, even though others map of the interval (as the logistic map) can be used for the same purpose. The associated dynamical system [24, 25] is defined by the equation on the interval [-1, 1]

$$x_{n+1} = 1 - a|x_n| \qquad (2.2)$$

Two tent maps are coupled in the following way, using a two dimensional coupling constant $\varepsilon = (\varepsilon_1, \varepsilon_2)$:

$$\begin{cases} x_{n+1} = (1-\varepsilon_1) f(x_n) + \varepsilon_1 f(y_n) \\ y_{n+1} = \varepsilon_2 f(x_n) + (1-\varepsilon_2) f(y_n) \end{cases} \qquad (2.3)$$

In this paper for the numerical studies we fix constant the ratio between $\varepsilon_1$ and $\varepsilon_2$. We chose to set it equal to 2. However, others ratios can also lead to good results and used as a multidimensional variable can be instrumental in the increasing of the number of dimensions of the systems.

$$\varepsilon_2 = 2\varepsilon_1 \qquad (2.4)$$

The coupling constant $\varepsilon$ varies from $(0, 0)$ to $(1, 1)$. When $\varepsilon = (0, 0)$ the maps are decoupled, when $\varepsilon = (1, 1)$ they are fully cross coupled. Generally authors do not consider very small values of $\varepsilon$ (as small as $10^{-7}$ for floating point numbers or $10^{-14}$ for double precision numbers), because it seems that the maps are quasi decoupled with those values. Hence no special effect of the coupling is expected.

In fact it is not the case and this very very small coupling constant allows the construction of very long periodic orbits, leading to sterling chaotic generators.

The dynamical system (2.3) can be described more generally by

$$\begin{pmatrix} x_{n+1} \\ y_{n+1} \end{pmatrix} = F\begin{pmatrix} x_n \\ y_n \end{pmatrix} = \begin{pmatrix} (1-\varepsilon_1) & \varepsilon_1 \\ \varepsilon_2 & (1-\varepsilon_2) \end{pmatrix} \cdot \begin{pmatrix} f(x_n) \\ f(y_n) \end{pmatrix} \qquad (2.5)$$

or

$$X_{n+1} = F(X_n) = A \cdot (\underline{f}(X_n)) \qquad (2.6)$$

with

$$X = \begin{pmatrix} x \\ y \end{pmatrix}, \qquad \underline{f}(X) = \begin{pmatrix} f(x) \\ f(y) \end{pmatrix} \quad \text{and} \quad A = \begin{pmatrix} (1-\varepsilon_1) & \varepsilon_1 \\ \varepsilon_2 & (1-\varepsilon_2) \end{pmatrix} \qquad (2.7)$$

$F$ is a map of the square $[-1, 1] \times [-1, 1]$ into itself.

### 2.1.2 Three and p-coupled symmetric map

In order to improve the length of the period and the convergence of the invariant measure towards a given measure, we will have to consider the coupling of three or more maps of the interval with the following coupling:

$$\begin{cases} x_{n+1} = (1-2\varepsilon_1)f(x_n) + \varepsilon_1 f(y_n) + \varepsilon_1 f(z_n) \\ y_{n+1} = \varepsilon_2 f(x_n) + (1-2\varepsilon_2)f(y_n) + \varepsilon_2 f(z_n) \\ z_{n+1} = \varepsilon_3 f(x_n) + \varepsilon_3 f(y_n) + (1-2\varepsilon_3)f(z_n) \end{cases} \quad (2.8)$$

and

$$\varepsilon_2 = 2\varepsilon_1, \quad \varepsilon_3 = 3\varepsilon_1 \quad (2.9)$$

(as mentioned above, others choices are possible). More generally the model of *p* coupled maps we introduce takes the form

$$\begin{pmatrix} x_{n+1}^1 \\ x_{n+1}^2 \\ \vdots \\ \vdots \\ x_{n+1}^p \end{pmatrix} = F \begin{pmatrix} x_n^1 \\ x_n^2 \\ \vdots \\ \vdots \\ x_n^p \end{pmatrix} = \begin{pmatrix} 1-(p-1)\varepsilon_1 & \varepsilon_1 & \cdots & \varepsilon_1 & \varepsilon_1 \\ \varepsilon_2 & 1-(p-1)\varepsilon_2 & \cdots & \varepsilon_2 & \varepsilon_2 \\ \vdots & & \ddots & \vdots & \vdots \\ \vdots & & & \ddots & \vdots \\ \varepsilon_p & \cdots & \cdots & \varepsilon_p & 1-(p-1)\varepsilon_p \end{pmatrix} \cdot \begin{pmatrix} f(x_n^1) \\ f(x_n^2) \\ \vdots \\ \vdots \\ f(x_n^p) \end{pmatrix} \quad (2.10)$$

with

$$\varepsilon_i = i\,\varepsilon_1 \qquad i = 2, \ldots, p \quad (2.11)$$

or

$$X_{n+1} = F(X_n) = A \cdot \left(\underline{f}(X_n)\right) \quad (2.12)$$

with

$$X = \begin{pmatrix} x^1 \\ \vdots \\ \vdots \\ x^p \end{pmatrix}, \quad \underline{f}(X) = \begin{pmatrix} f(x^1) \\ \vdots \\ \vdots \\ f(x^p) \end{pmatrix} \quad (2.13)$$

and

$$A = \begin{pmatrix} 1-(p-1)\varepsilon_1 & \varepsilon_1 & \cdots & \varepsilon_1 & \varepsilon_1 \\ \varepsilon_2 & 1-(p-1)\varepsilon_2 & \cdots & \varepsilon_2 & \varepsilon_2 \\ \vdots & & \ddots & \vdots & \vdots \\ \vdots & & & \ddots & \vdots \\ \varepsilon_p & \cdots & \cdots & \varepsilon_p & 1-(p-1)\varepsilon_p \end{pmatrix} \quad (2.14)$$

In this case *F* is a map of $[-1, 1]^p$ into itself.

### 2.2 Numerical results for chaotic generators

In this subsection we verify for 2, 3 and 4-coupled symmetric tent maps that the chaotic numbers produced by means of the equation (2.12) are equally distributed over the interval [-1, 1].

That is, the invariant measure of each component converges towards the Lebesgue measure. We conduct numerical experiments testing first the influence of the number of iterations and the value of the discretization of the approximation of the invariant measure. We test next in subsection 2.2.3 the influence of the initial conditions, and finally in subsection 2.2.4 we show that the series of chaotic numbers produced by each component is independent of the others.

### 2.2.1 Approximated distribution function

In order to compute numerically an approximation of the invariant measure [21] also called the probability distribution function $P_N(x)$ linked to the one dimensional map $f$ we build a regular partition of $M$ small intervals (boxes) of the considered interval [-1, 1] :

$$r_i = [s_i, s_{i+1}[, \quad i = 0, M-2 \quad (2.15)$$

$$r_{M-1} = [s_{M-1}, 1] \quad (2.16)$$

$$s_i = -1 + \frac{2i}{M} \quad i = 0, M \quad (2.17)$$

the length of which is:

$$s_{i+1} - s_i = \frac{2}{M} \quad (2.18)$$

We collect all iterates $f^{(n)}(x)$ belonging to these boxes (after a transient regime of $q$ iterations decided a priori, *i.e.* the first $q$ iterates are neglected). Once the computation of $N+q$ iterates is completed, the relative number of iterates with respect to $N/M$ in each box $r_i$ represents the value $P_N(s_i)$. The approximated probability distribution function $P_N(x)$ defined in this article is then a step function, with $M$ steps. As $M$ can vary in the next sections, we define :

$$P_{M,N}(s_i) = \frac{1}{2}\frac{M}{N}(\#r_i) \quad (2.19)$$

where $\#r_i$ is the number of iterates belonging to the interval $r_i$ and the constant $\frac{1}{2}$ allows the normalisation of $P_{M,N}(x)$ on the interval [-1, 1].

$$P_{M,N}(x) = P_{M,N}(s_i) \quad \forall x \in r_i \quad (2.20)$$

In the case of coupled maps, we are interested by the distribution of the component $x^1$, …, $x^p$ of $X$ rather than the distribution of the variable $X$ itself in $[-1, 1]^p$. Then we will consider the approximated probability distribution function $P_N(x^j)$ associated to one within several components of $F(X)$ defined by (2.12) which are one-dimensional maps.

### 2.2.2 Distribution of iterates of 2, 3 and 4-coupled symmetric tent maps

In this section, the numerical experiments are performed on several computers involving different microprocessors of Advanced Micro Devices (AMD) and Intel (Centrino and dual core) technologies in order to check the portability of the algorithms we propose. In the same goal the package is written using many versions of Borland C.

Double precision numbers are used (see [21] for tests with floating points). We fix $\varepsilon = 10^{-14}$.

As intuitively expected, the density of iterates of each component of 2-coupled and *p*-coupled symmetric tent map converges towards the Lebesgue measure when $\varepsilon_1$ converges towards 0. We call discretization number ($N_{\text{disc}}$) the number *M* of boxes of the interval and number of iterates ($N_{\text{iter}}$) the number *N*.

The asymptotic properties of dynamical systems intuitively imply that for a fixed value of $N_{\text{disc}}$ when the number $N_{\text{iter}}$ increases, the discrepancy between $P_{N_{disc}, N_{iter}}(x)$ and the Lebesgue measure converge towards 0, except if there exist one or many periodic orbits of finite length lower than $N_{\text{iter}}$ which capture the iterates. In this case whatsoever the value of $N_{\text{iter}}$ is, the approximated distribution function converges to the distribution function of the periodic orbit if it is unique or to some average of the distribution functions of the periodic orbits observed if not.

The discrepancies $E_1$ (i.e. in norm $L_1$) and $E_2$ (in norm $L_2$) between $P_{N_{disc}, N_{iter}}(x)$ and the Lebesgue measure are defined by:

$$E_1(N_{disc}, N_{iter}) = \left\| P_{N_{disc}, N_{iter}}(x) - 0.5 \right\|_{L_1} \qquad (2.21)$$

$$E_2(N_{disc}, N_{iter}) = \left\| P_{N_{disc}, N_{iter}}(x) - 0.5 \right\|_{L_2} \qquad (2.22)$$

Figures 1 and 2 show the errors $E_1(N_{disc}, N_{iter})$ and $\left(E_2(N_{disc}, N_{iter})\right)^2$ versus the number of iterates of the approximated distribution functions with respect to the first variable $x^1$ for 2, 3 and 4-coupled symmetric tent map. $N_{\text{disc}}$ is fixed to $10^4$, $\varepsilon_1$ to $10^{-14}$, $N_{\text{iter}}$ varies from $10^5$ to $10^{11}$ for the 2-coupled case and to $3.10^{12}$ for the 3 and 4-coupled one. The corresponding numerical results are displayed in Table 1 $E_1(N_{disc}, N_{iter})$ for and Table 2 for $\left(E_2(N_{disc}, N_{iter})\right)^2$.

**Remark:** in order to made easier the comparison of the results, we display the square of the discrepancy $E_2$ instead of $E_2$ itself, the discrepancy being divided by 10 each time the number of iterations is multiplied by 10.

Equivalent results are obtained for the variables $x^2$, $x^3$ or $x^4$.

No periodic solutions are observed up to $3.10^{12}$ iterates (even up to $10^{13}$ iterates as tested in section 3). This point is very important for producing chaotic numbers, because the use of a computer discretizes the phase space of a dynamical system, canceling (at least) its asymptotic properties. Every orbit is periodic according to the finite number of states (*i. e.,* the number of double precision numbers belonging to the *p*-dimensional interval $[-1, 1]^p$). However, if the period of the realized sequence is long enough, these properties reasonably survive as a chaotic transient.

One can observe that for 3 and 4-coupled equations the convergence is excellent up to $3.10^{12}$ iterates. For 2-coupled equations the convergence seems lower bounded by a minimal error.

There is no significant difference between 3 and 4-coupled equations, the numerical experiments have to be pursued up to $10^{13}$ or $10^{14}$ in order to discriminate the results.

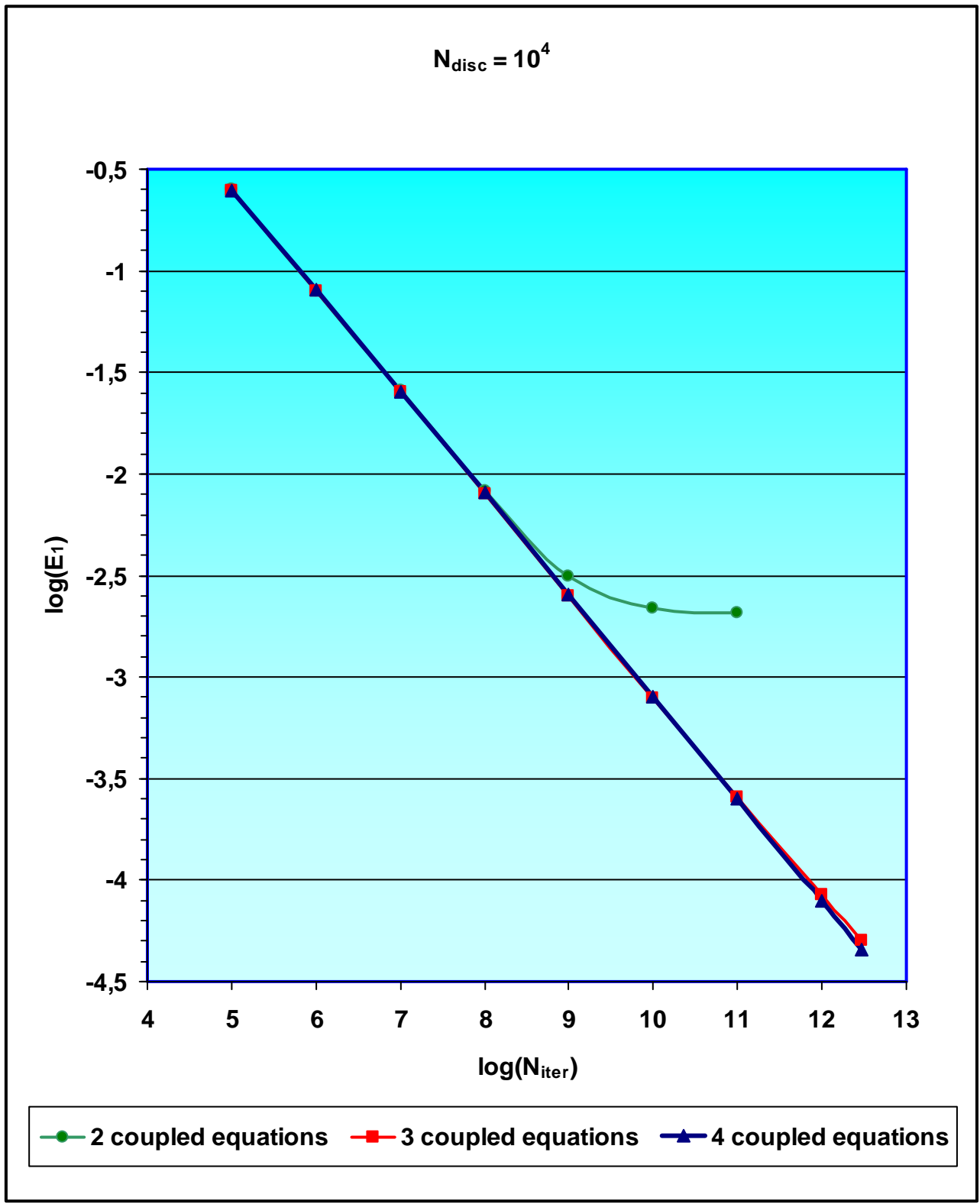

**Figure 1.** Error $E_1$ for 2, 3 and 4-coupled Symmetric Tent Maps. Computations done using double precision numbers (~14-15 digits), $\varepsilon_i = i.\varepsilon_1$, $\varepsilon_1 = 10^{-14}$, $N_{disc} = 10^4$.

Initial values $x^1_0 = 0.330000013113$, $x^2_0 = 0.338756413113$, $x^3_0 = 0.331353442113$, $x^4_0 = 0.333213583113$.

| $N_{iter}$ | $E_1(N_{disc}, N_{iter})$ **2-coupled equation** | $E_1(N_{disc}, N_{iter})$ **3-coupled equation** | $E_1(N_{disc}, N_{iter})$ **4-coupled equation** |
|---|---|---|---|
| $10^5$ | 0.25071335 | 0.25035328 | 0.2499133 |
| $10^6$ | 0.079655103 | 0.079437105 | 0.080739109 |
| $10^7$ | 0.025794703 | 0.025343302 | 0.025266304 |
| $10^8$ | 0.0081966502 | 0.0079505501 | 0.0080771501 |
| $10^9$ | 0.003147609 | 0.002513533 | 0.002562893 |
| $10^{10}$ | 0.0021717461 | 0.0007908719 | 0.0007970199 |
| $10^{11}$ | 0.0020550967 | 0.00025791013 | 0.00025241399 |
| $10^{12}$ | | $8.4195287 \cdot 10^{-5}$ | $7.8803383 \cdot 10^{-5}$ |
| $3 \cdot 10^{12}$ | | $5.0625114 \cdot 10^{-5}$ | $4.5317128 \cdot 10^{-5}$ |

**Table 1.** Error $E_1$ for 2, 3 and 4-coupled Symmetric Tent Maps. Computations done using double precision numbers (~14-15 digits), $\varepsilon_i = i \cdot \varepsilon_1$, $\varepsilon_1 = 10^{-14}$, $N_{disc} = 10^4$.
Initial values $x^1_0 = 0.330000013113$, $x^2_0 = 0.338756413113$, $x^3_0 = 0.331353442113$, $x^4_0 = 0.333213583113$.

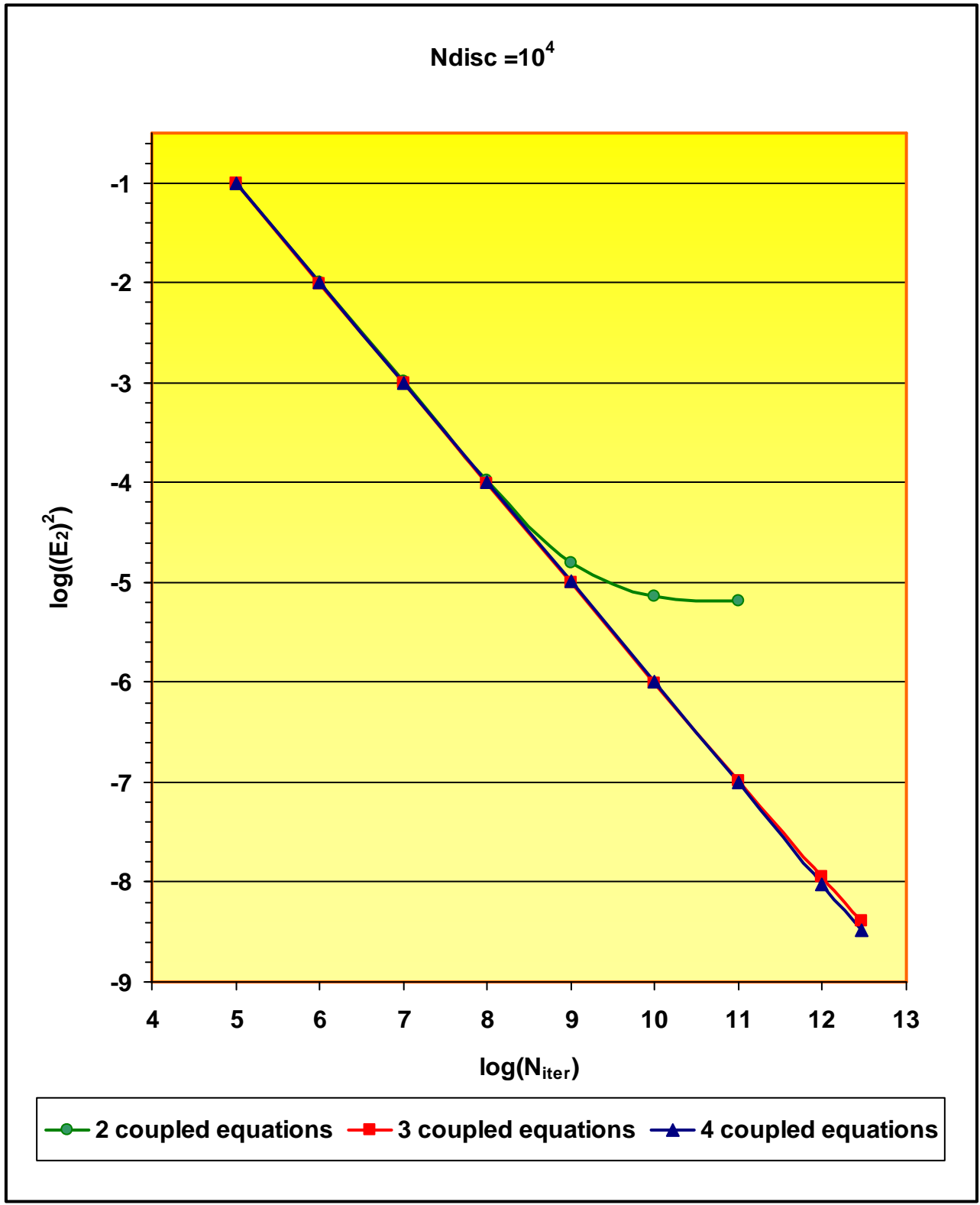

**Figure 2.** Error $E_2^2$ for 2, 3 and 4-coupled Symmetric Tent Maps. Computations done using double precision numbers (~14-15 digits), $\varepsilon_i = i.\varepsilon_1$, $\varepsilon_1 = 10^{-14}$, $N_{\text{disc}} = 10^4$.

Initial values $x^1_0 = 0.330000013113$, $x^2_0 = 0.338756413113$, $x^3_0 = 0.331353442113$, $x^4_0 = 0.333213583113$.

| $N_{iter}$ | $E_2^2(N_{disc}, N_{iter})$ 2-coupled equation | $E_2^2(N_{disc}, N_{iter})$ 3-coupled equation | $E_2^2(N_{disc}, N_{iter})$ 4-coupled equation |
|---|---|---|---|
| $10^5$ | 0.100199 | 0.099820996 | 0.099610992 |
| $10^6$ | 0.01006199 | 0.0098781898 | 0.01022057 |
| $10^7$ | 0.0010442081 | 0.0010014581 | 0.0010055967 |
| $10^8$ | 0.0001055816 | $9.8853067 \cdot 10^{-5}$ | 0.00010197872 |
| $10^9$ | $1.567597 \cdot 10^{-5}$ | $1.0047459 \cdot 10^{-5}$ | $1.0326474 \cdot 10^{-5}$ |
| $10^{10}$ | $7.3577797 \cdot 10^{-6}$ | $9.7251536 \cdot 10^{-7}$ | $9.9932242 \cdot 10^{-7}$ |
| $10^{11}$ | $6.6338453 \cdot 10^{-6}$ | $1.0434293 \cdot 10^{-7}$ | $1.0070523 \cdot 10^{-7}$ |
| $10^{12}$ | | $1.116009 \cdot 10^{-8}$ | $9.6166733 \cdot 10^{-9}$ |
| $3 \cdot 10^{12}$ | | $4.0443118 \cdot 10^{-9}$ | $3.2530773 \cdot 10^{-9}$ |

**Table 2.** Error $E_2^2$ for 2, 3 and 4-coupled Symmetric Tent Maps. Computations done using double precision numbers (~14-15 digits), $\varepsilon_i = i.\varepsilon_1$, $\varepsilon_1 = 10^{-14}$, $N_{disc} = 10^4$.

Initial values $x^1_0 = 0.330000013113$, $x^2_0 = 0.338756413113$, $x^3_0 = 0.331353442113$, $x^4_0 = 0.333213583113$.

In Figure 3, 4 and 5 we display the mutual influence of both $N_{iteration}$ and $N_{discretization}$ on the errors in $L_1$ and $L_2$ norm. The results show a tremendous regularity.

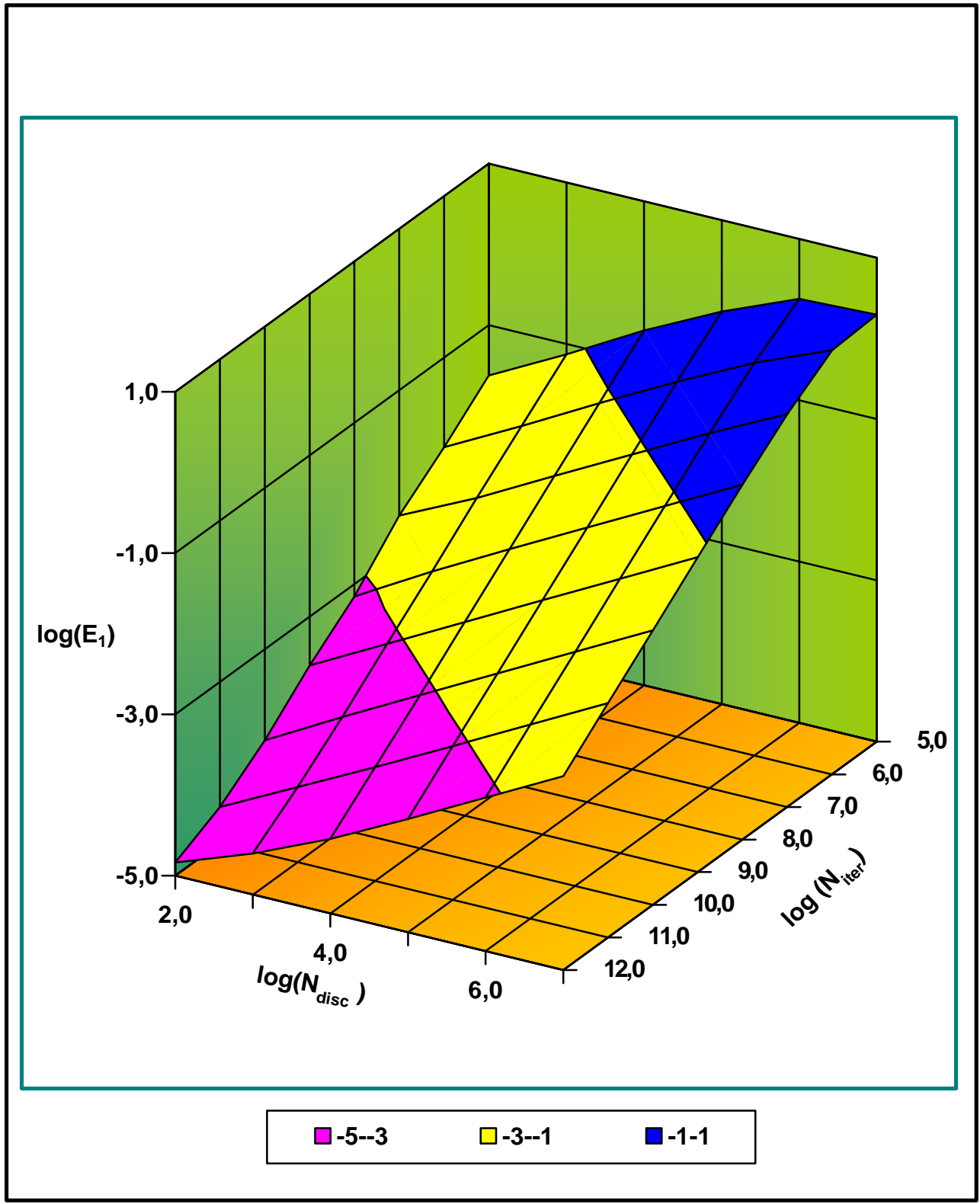

**Figure 3.** Error $E_1$ for 3-coupled Symmetric Tent Maps. Computations done using double precision numbers (~14-15 digits) with respect to both $N_{iteration}$ and $N_{discretization}$, $\varepsilon_i = i.\varepsilon_1$, $\varepsilon_1 = 10^{-14}$, $N_{iter} = 10^5$ to $10^{12}$, $N_{disc} = 10^2$ to $10^7$.

Initial values $x^1_0 = 0.330000013113$, $x^2_0 = 0.338756413113$, $x^3_0 = 0.331353442113$.

| $N_{iter}$ \ $N_{disc}$ | $E_1(N_{disc}, N_{iter})$ $10^2$ | $E_1(N_{disc}, N_{iter})$ $10^3$ | $E_1(N_{disc}, N_{iter})$ $10^4$ |
|---|---|---|---|
| $10^5$ | 0.023590236 | 0.074390944 | 0.25035328 |
| $10^6$ | 0.0077829878 | 0.024115036 | 0.079437105 |
| $10^7$ | 0.0027963003 | 0.0078734998 | 0.025343302 |
| $10^8$ | 0.00070102901 | 0.0024396098 | 0.0079505501 |
| $10^9$ | 0.00024907298 | 0.00078846501 | 0.002513533 |
| $10^{10}$ | $7.4041294 \cdot 10^{-5}$ | 0.0002472693 | 0.0007908719 |
| $10^{11}$ | $2.821469 \cdot 10^{-5}$ | $8.540793 \cdot 10^{-5}$ | 0.00025791013 |
| $10^{12}$ | $1.4600127 \cdot 10^{-5}$ | $3.2358931 \cdot 10^{-5}$ | $8.4195287 \cdot 10^{-5}$ |

| $N_{iter}$ \ $N_{disc}$ | $E_1(N_{disc}, N_{iter})$ $10^5$ | $E_1(N_{disc}, N_{iter})$ $10^6$ | $E_1(N_{disc}, N_{iter})$ $10^7$ |
|---|---|---|---|
| $10^5$ | 0.73832 | 1.810124 | 1.9801114 |
| $10^6$ | 0.24974733 | 0.735708 | 1.809666 |
| $10^7$ | 0.079959311 | 0.25029673 | 0.7353684 |
| $10^8$ | 0.02518029 | 0.079508971 | 0.25000429 |
| $10^9$ | 0.008005619 | 0.025207567 | 0.079757051 |
| $10^{10}$ | 0.0025136649 | 0.0079736449 | 0.025230797 |
| $10^{11}$ | 0.00080110625 | 0.002522144 | 0.0079771447 |
| $10^{12}$ | 0.00025407246 | 0.00079907514 | 0.0025234708 |

**Table 3.** Error $E_1$ for 3-coupled Symmetric Tent Maps. Computations done using double precision numbers (~14-15 digits) with respect to both $N_{iteration}$ and $N_{discretization}$, $\varepsilon_i = i.\varepsilon_1$, $\varepsilon_1 = 10^{-14}$, $N_{iter} = 10^5$ to $10^{12}$, $N_{disc} = 10^2$ to $10^7$.
Initial values $x^1_0 = 0.330000013113$, $x^2_0 = 0.338756413113$, $x^3_0 = 0.331353442113$.

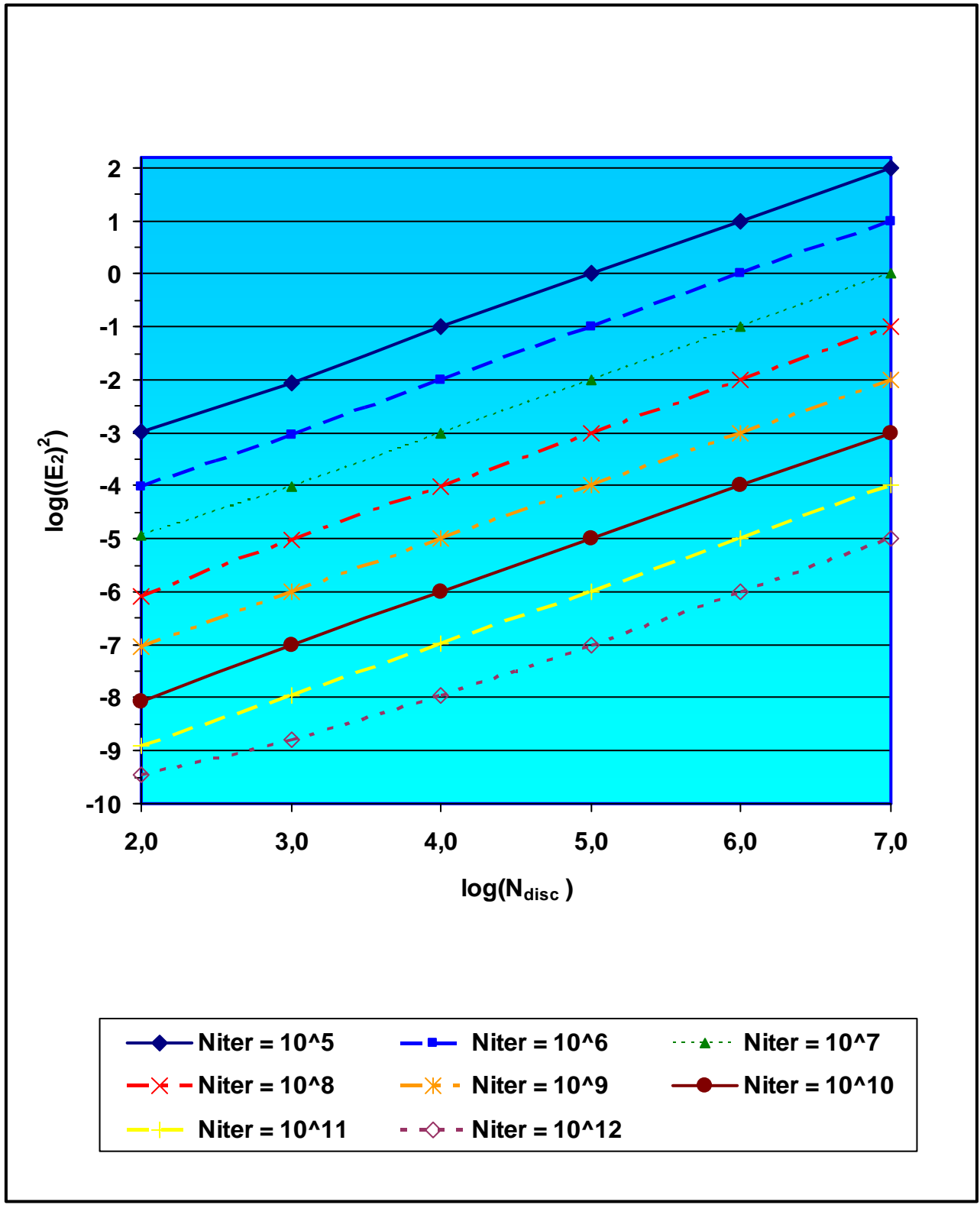

**Figure 4.** Error $E_2^2$ for 3-coupled Symmetric Tent Maps. Computations done using double precision numbers (~14-15 digits) with respect to both $N_{\text{iteration}}$ and $N_{\text{discretization}}$, $\varepsilon_i = i.\varepsilon_1$, $\varepsilon_1 = 10^{-14}$, $N_{\text{iter}} = 10^5$ to $10^{12}$, $N_{\text{disc}} = 10^2$ to $10^7$.
Initial values $x^1_0 = 0.330000013113$, $x^2_0 = 0.338756413113$, $x^3_0 = 0.331353442113$.

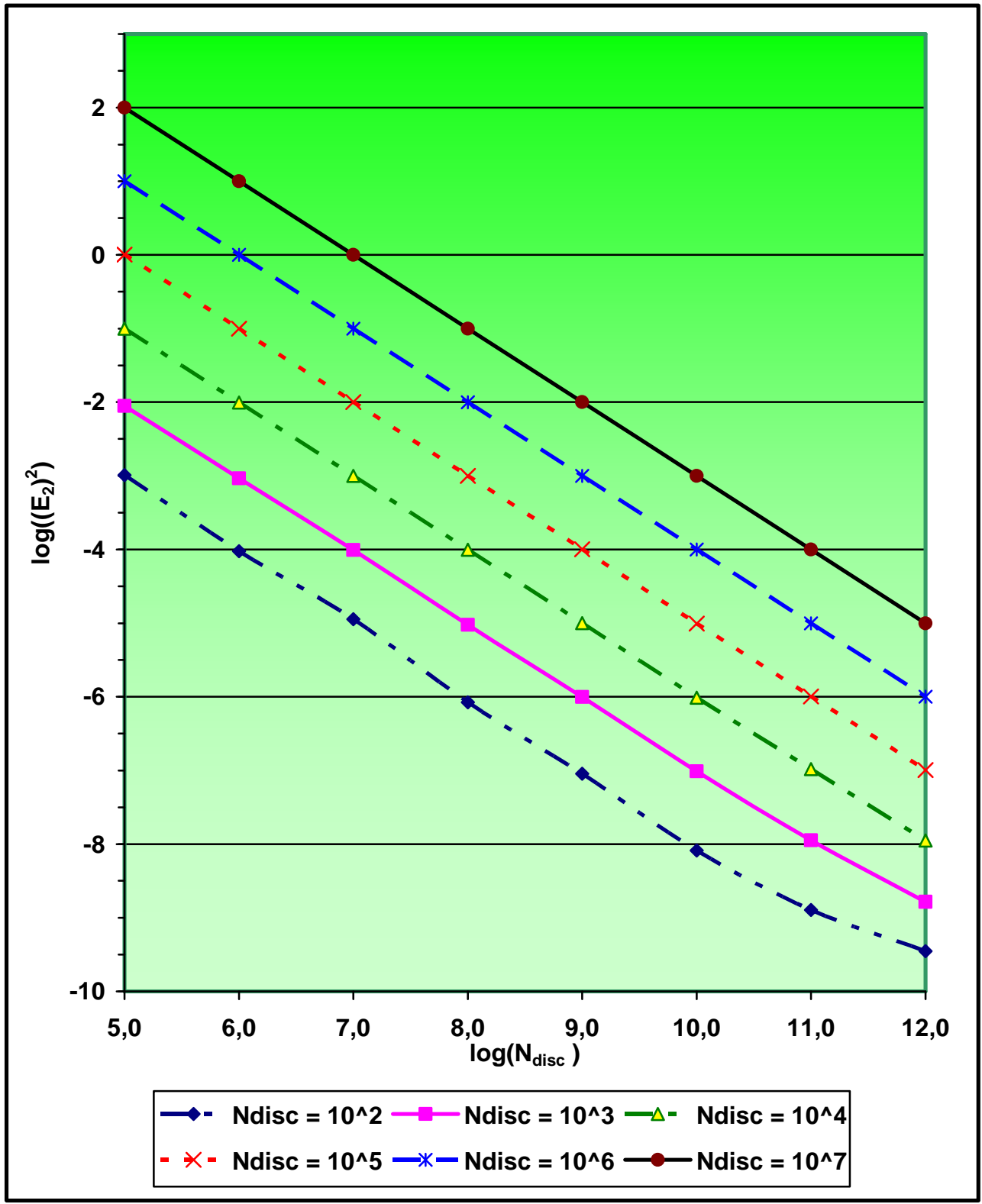

**Figure 5.** Error $E_2^2$ for 3-coupled Symmetric Tent Maps. Computations done using double precision numbers (~14-15 digits) with respect to both $N_{\text{iteration}}$ and $N_{\text{discretization}}$, $\varepsilon_i = i.\varepsilon_1$, $\varepsilon_1 = 10^{-14}$, $N_{\text{iter}} = 10^5$ to $10^{12}$, $N_{\text{disc}} = 10^2$ to $10^7$.

Initial values $x^1_0 = 0.330000013113$, $x^2_0 = 0.338756413113$, $x^3_0 = 0.331353442113$.

| $N_{iter}$ \ $N_{disc}$ | $E_2^2(N_{disc}, N_{iter})$ $10^2$ | $E_2^2(N_{disc}, N_{iter})$ $10^3$ | $E_2^2(N_{disc}, N_{iter})$ $10^4$ |
|---|---|---|---|
| $10^5$ | 0.0010228504 | 0.0088666772 | 0.099820996 |
| $10^6$ | $9.4971289 \cdot 10^{-5}$ | 0.00092621485 | 0.0098781898 |
| $10^7$ | $1.1384949 \cdot 10^{-5}$ | $9.827981 \cdot 10^{-5}$ | 0.0010014581 |
| $10^8$ | $8.3729789 \cdot 10^{-7}$ | $9.5104663 \cdot 10^{-6}$ | $9.8853067 \cdot 10^{-5}$ |
| $10^9$ | $9.0592196 \cdot 10^{-8}$ | $9.9378865 \cdot 10^{-7}$ | $1.0047459 \cdot 10^{-5}$ |
| $10^{10}$ | $8.1956005 \cdot 10^{-9}$ | $9.7587086 \cdot 10^{-8}$ | $9.7251536 \cdot 10^{-7}$ |
| $10^{11}$ | $1.275616 \cdot 10^{-9}$ | $1.1261091 \cdot 10^{-8}$ | $1.0434293 \cdot 10^{-7}$ |
| $10^{12}$ | $3.521566 \cdot 10^{-10}$ | $1.6424142 \cdot 10^{-9}$ | $1.116009 \cdot 10^{-8}$ |

| $N_{iter}$ \ $N_{disc}$ | $E_2^2(N_{disc}, N_{iter})$ $10^5$ | $E_2^2(N_{disc}, N_{iter})$ $10^6$ | $E_2^2(N_{disc}, N_{iter})$ $10^7$ |
|---|---|---|---|
| $10^5$ | 1.0118502 | 10.047321 | 100.12102 |
| $10^6$ | 0.0998531 | 1.000171 | 10.000052 |
| $10^7$ | 0.010018681 | 0.09999495 | 0.9987099 |
| $10^8$ | 0.00099451547 | 0.0099463249 | 0.099816961 |
| $10^9$ | 0.00010078117 | 0.00099783362 | 0.010007831 |
| $10^{10}$ | $9.923092 \cdot 10^{-6}$ | $9.9847198 \cdot 10^{-5}$ | 0.00099995901 |
| $10^{11}$ | $1.007156 \cdot 10^{-6}$ | $9.9997448 \cdot 10^{-6}$ | $9.9965829 \cdot 10^{-5}$ |
| $10^{12}$ | $1.0137626 \cdot 10^{-7}$ | $1.0018743 \cdot 10^{-6}$ | $1.0000158 \cdot 10^{-5}$ |

**Table 4.** Error $E_2^2$ for 3-coupled Symmetric Tent Maps. Computations done using double precision numbers (~14-15 digits) with respect to both $N_{\text{iteration}}$ and $N_{\text{discretization}}$, $\varepsilon_i = i.\varepsilon_1$, $\varepsilon_1 = 10^{-14}$, $N_{\text{iter}} = 10^5$ to $10^{12}$, $N_{\text{disc}} = 10^2$ to $10^7$.

Initial values $x^1_0 = 0.330000013113$, $x^2_0 = 0.338756413113$, $x^3_0 = 0.331353442113$.

### 2.2.3 Impact of the initial values on the results

It is well known that the choice of the seed of a PRNG is very important. Some seed can lead to the collapse of the period of the computed random numbers. In order to check if the choice of the initial condition (equivalent to the choice of the seed of a PRNG) is dramatically for the previous results, we have tested two series of different initial values.

Figure 6 shows the distribution of the error $E_1$ for 500,000 initial values for 4-coupled symmetric tent maps. The computations are done using double precision numbers (~14-15 digits), $\varepsilon_i = i.\varepsilon_1$, $\varepsilon_1 = 10^{-14}$, $N_{\text{iter}} = 10^6$, $N_{\text{disc}} = 10^2$.

The initial values are selected following:

$$x^1_{0,k} = -0.92712 + 10^{-7} \times k, \quad x^2_{0,k} = -0.9183636 + 10^{-7} \times 7k,$$
$$x^3_{0,k} = -0.92576657 + 10^{-7} \times 13k, \quad x^4_{0,k} = -0.92390643 + 10^{-7} \times 17k,$$
$$k = 1 \text{ to } 500,000.$$

The distribution follows more or less a Gaussian distribution, maximal and minimal results are displayed in Table 5.

Others results for 3-coupled equations for the sequence of initial values:

$$x^1_{0,k} = -0.92712 + 10^{-6} \times k, \quad x^2_{0,k} = -0.9183636 + 10^{-6} \times 7k, \quad x^3_{0,k} = -0.92576657 + 10^{-6} \times 13k.$$
$$k = 1 \text{ to } 100,000$$

are plotted in Figure 7 for $N_{\text{disc}} = 10^2$ to $10^4$.

All these results confirm that the families of chaotic attractor we have introduced are robust versus the choice of the initial seed.

| $N_{\text{disc}}$ | $10^2$ | $10^3$ | $10^4$ |
|---|---|---|---|
| $\min E_1(N_{disc}, N_{iter})$ | 0.004002104 | 0.020740039 | 0.075152102 |
| $\max E_1(N_{disc}, N_{iter})$ | 0.013871994 | 0.030116008 | 0.084384101 |
| $\min E_2^2(N_{disc}, N_{iter})$ | 0.000027527656 | 0.00067688535 | 0.0089216578 |
| $\max E_2^2(N_{disc}, N_{iter})$ | 0.00028341257 | 0.0014350609 | 0.011071942 |

**Table 5.** Minimal and maximal values of the $E_1$ and $E_2^2$ errors for 500,000 initial values for 4-coupled symmetric tent maps. Computations done using double precision numbers (~14-15 digits), $\varepsilon_i = i.\varepsilon_1$, $\varepsilon_1 = 10^{-14}$, $N_{\text{iter}} = 10^6$, $N_{\text{disc}} = 10^2$ to $10^4$.
Initial values $x^1_{0,k} = -0.92712 + 10^{-7} \times k$, $x^2_{0,k} = -0.9183636 + 10^{-7} \times 7k$, $x^3_{0,k} = -0.92576657 + 10^{-7} \times 13k$, $x^4_{0,k} = -0.92390643 + 10^{-7} \times 17k$.

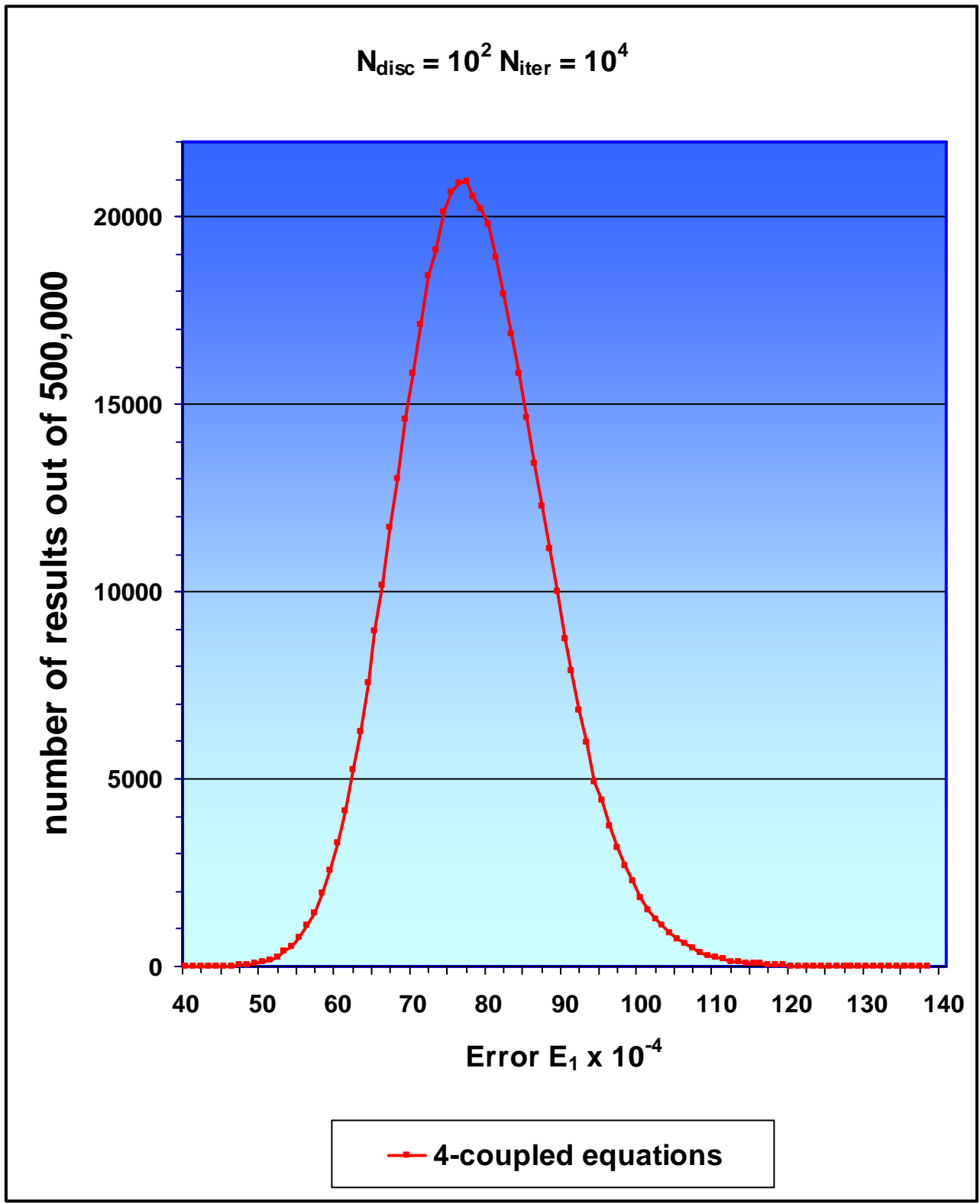

**Figure 6.** Distribution of the error $E_1$ for 500.000 initial values for 4-coupled symmetric tent maps. Computations done using double precision numbers (~14-15 digits), $\varepsilon_i = i.\varepsilon_1$, $\varepsilon_1 = 10^{-14}$, $N_{\text{iter}} = 10^6$, $N_{\text{disc}} = 10^2$.
Initial values $x^1_{0,k} = -0.92712 + 10^{-7} \times k$, $x^2_{0,k} = -0.9183636 + 10^{-7} \times 7k$, $x^3_{0,k} = -0.92576657 + 10^{-7} \times 13k$, $x^4_{0,k} = -0.92390643 + 10^{-7} \times 17k$.

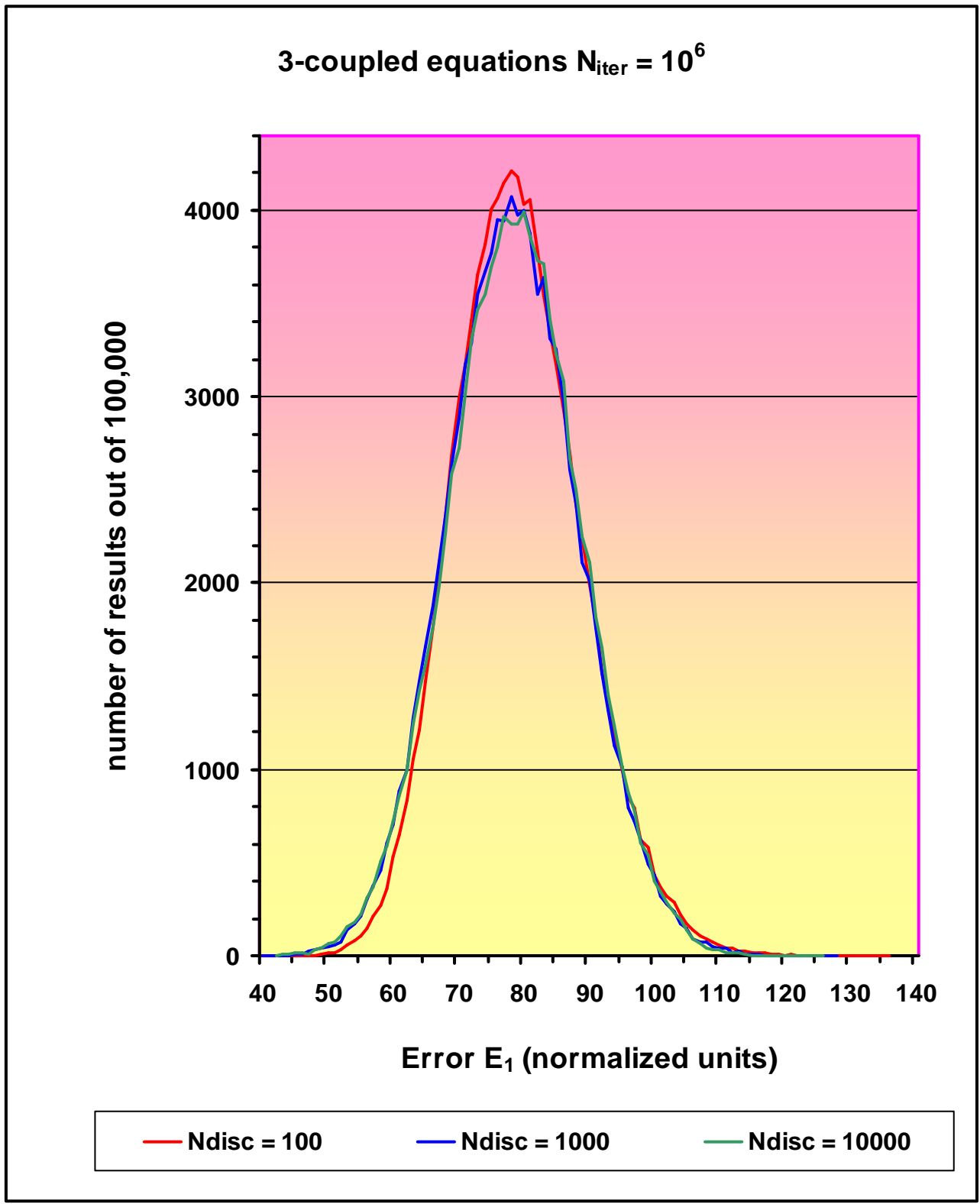

**Figure 7.** Distribution of the error $E_1$ for 100.000 initial values for 3-coupled symmetric tent maps. Computations done using double precision numbers (~14-15 digits), $\varepsilon_i = i.\varepsilon_1$, $\varepsilon_1 = 10^{-14}$, $N_{iter} = 10^6$, $N_{disc} = 10^2$ to $10^4$.
Initial values $x^1_{0,k} = -0.92712 + 10^{-6} \times k$, $x^2_{0,k} = -0.9183636 + 10^{-6} \times 7k$, $x^3_{0,k} = -0.92576657 + 10^{-6} \times 13k$.

### 2.2.4 Independency of the chaotic subsequences generated by each component

In Section 3 we will propose chaotic samplings of the chaotic sequences generated by (2.10) to enhance the properties of this chaotic number generator. The key feature of these enhanced chaotic number generators being their use of chaotic numbers themselves in order to do the sampling process. The main idea leading to this particular sampling is that the series of chaotic numbers produced by each component is independent of the others.

We need before to verify this independency.

We recall that $X_n = \begin{pmatrix} x_n^1 \\ \vdots \\ \vdots \\ x_n^p \end{pmatrix}$. Let consider now the coordinates of the iterated points $X_0, X_1, X_2, \cdots, X_n, X_{n+1}, \cdots$ of the multidimensional map $F$ defined by (2.12). In order to control that they are uncorrelated, we plot every pair of coordinates of this sequence in the phase subspace $(x^l, x^m)$ imbedded in the phase space $[-1, 1]^p$ and we check if they are uniformly distributed in the square $[-1, 1]^2$.

If no particular pattern is displayed and if the difference between the distribution of these points later called the correlation distribution function $C_N(x, y)$ converges towards the uniform distribution on the square when the number of iterations goes to the infinity, we can conclude the independency or the uncorrelation of the sequences of numbers generated by each component of the iterated points.

In order to compute numerically an approximation of the correlation distribution function $C_N(x, y)$ we build a regular partition of $M^2$ small squares (boxes) of the square $[-1, 1]^2$ imbedded in the phase subspace $(x^l, x^m)$ :

$$r_{i,j} = [s_i, s_{i+1}[ \times [t_j, t_{j+1}[ \quad , \; i,j = 0, M-2 \quad (2.23)$$

$$r_{M-1,j} = [s_{M-1}, 1] \times [t_j, t_{j+1}[ \quad , \; j = 0, M-2 \quad (2.24)$$

$$r_{i,M-1} = [s_i, s_{i+1}[ \times [t_{M-1}, 1] \quad , \; i = 0, M-2 \quad (2.25)$$

$$r_{M-,M-1} = [s_{M-1}, 1] \times [t_{M-1}, 1] \quad (2.26)$$

$$s_i = -1 + \frac{2i}{M} \quad i = 0, M \quad (2.27)$$

$$t_j = -1 + \frac{2j}{M} \quad j = 0, M \quad (2.28)$$

the measure of the area of each box is :

$$(s_{i+1} - s_i) \cdot (t_{i+1} - t_i) = \left(\frac{2}{M}\right)^2 \quad (2.29)$$

We collect all iterated points $(x_n^l, x_n^m)$ belonging to these boxes (after a transient regime of $q$ iterations decided a priori, *i.e.* the first $q$ iterates are neglected). Once the computation of $N + q$ iterates is completed, the relative number of iterates with respect to $N/M^2$ in each box $r_{i,j}$ represents the value $C_N(s_i, t_j)$. The approximated probability distribution function $C_N(x, y)$ defined in this article is then a 2-dimensional step function, with $M^2$ steps. As $M$ can vary in the next sections, we define :

$$C_{M,N}(s_i, t_j) = \frac{1}{4} \frac{M^2}{N} (\# r_{i,j}) \quad (2.30)$$

where $\# r_{i,j}$ is the number of iterates belonging to the square $r_{i,j}$ and the constant $\frac{1}{4}$ allows the normalisation of $C_{M,N}(x, y)$ on the square $[-1, 1]^2$.

$$C_{M,N}(x, y) = C_{M,N}(s_i, t_j) \quad \forall (x, y) \in r_{i,j} \quad (2.31)$$

The discrepancies $E_{C1}$ (i.e. in norm $L_1$) and $E_{C2}$ (in norm $L_2$) between $C_{N_{disc}, N_{iter}}(x, y)$ and the uniform distribution on the square are defined by:

$$E_{C1}(N_{disc}, N_{iter}) = \left\| C_{N_{disc}, N_{iter}}(x, y) - 0.25 \right\|_{L_1} \quad (2.32)$$

$$E_{C2}(N_{disc}, N_{iter}) = \left\| C_{N_{disc}, N_{iter}}(x, y) - 0.25 \right\|_{L_2} \quad (2.33)$$

Figure 8 displays the error $E_{C1}(N_{disc}, N_{iter})$ versus the number of iterated points of the approximated correlation function between the first and the second components $(x^1, x^2)$ for the 4-coupled symmetric tent map. $N_{disc}$ is fixed to $10^2 \times 10^2$, $\varepsilon_1$ to $10^{-14}$, $N_{iter}$ varies from $10^5$ to $10^{11}$. The corresponding numerical results are displayed in Table 6.

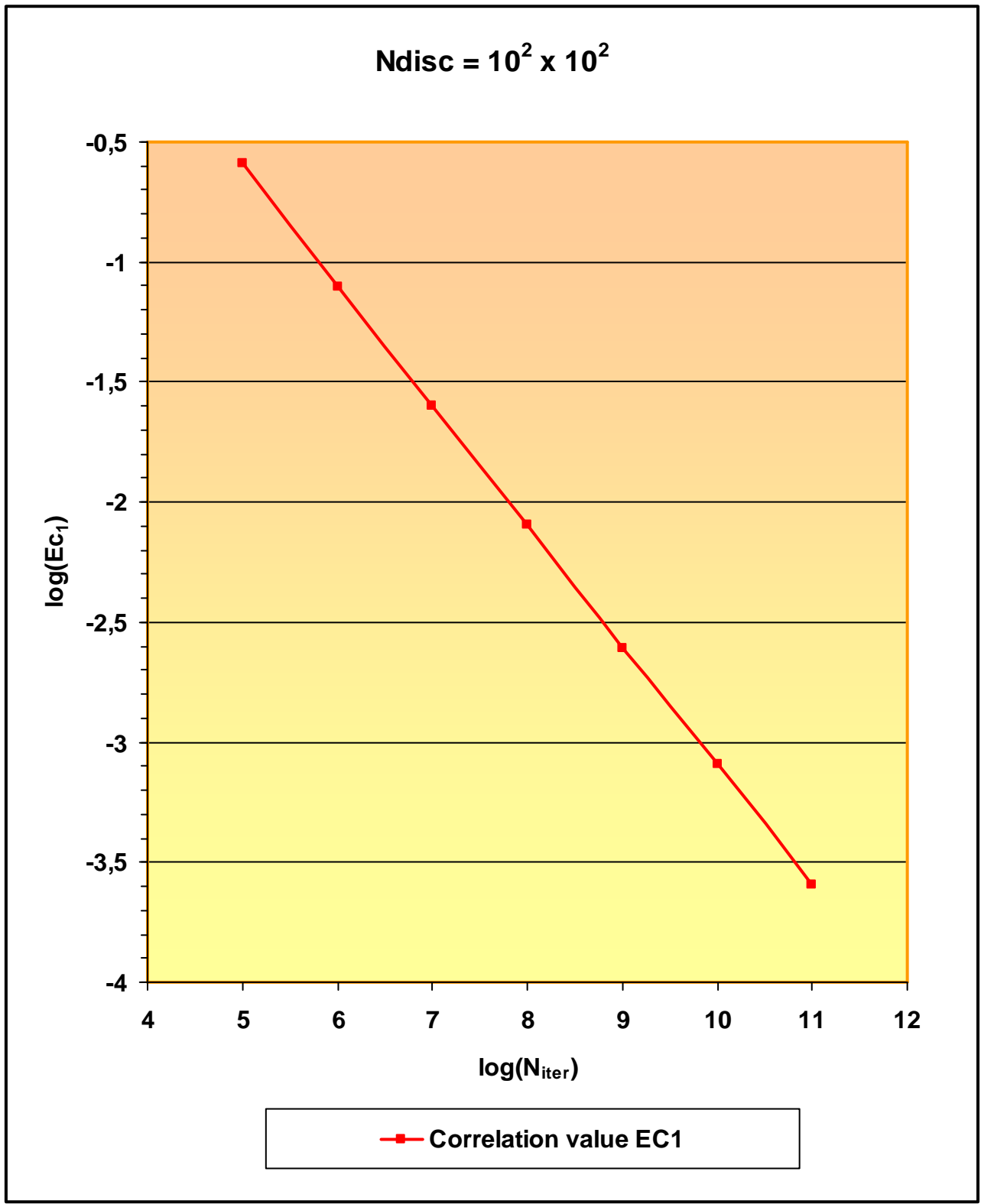

**Figure 8.** Error $E_{C_1}$ for the first and the second components $(x^1, x^2)$ of the 4-coupled symmetric tent map. $N_{disc}$ is fixed to $10^2 \times 10^2$, $\varepsilon_i = i.\varepsilon_1$, $\varepsilon_1 = 10^{-14}$, $N_{iter}$ varies from $10^5$ to $10^{11}$. Computations are done using double precision numbers (~14-15 digits).
Initial values $x^1_0 = 0.330$, $x^2_0 = 0.3387564$, $x^3_0 = 0.50492331$, $x^4_0 = 0.0$.

| $N_{\text{iter}}$ | $E_{C1(N_{disc},N_{iter})}(x^1, x^2)$ |
|---|---|
| $10^5$ | 0.2573333 |
| $10^6$ | 0.078763097 |
| $10^7$ | 0.025002305 |
| $10^8$ | 0.0080488902 |
| $10^9$ | 0.002477239 |
| $10^{10}$ | 0.0008041077 |
| $10^{11}$ | 0.00025640111 |

**Table 6.** Error $E_{C1}$ for the first and the second components $(x^1, x^2)$ of the 4-coupled symmetric tent map. $N_{\text{disc}} = 10^2 \times 10^2$, $\varepsilon_i = i.\varepsilon_1$, $\varepsilon_1 = 10^{-14}$, $N_{\text{iter}}$ varies from $10^5$ to $10^{11}$. Computations are done using double precision numbers (~14-15 digits).

Initial values $x^1_0 = 0.330$, $x^2_0 = 0.3387564$, $x^3_0 = 0.50492331$, $x^4_0 = 0.0$.

The difference between the correlation distribution function $C_{\text{NDISC}}(x^1, x^3)$ and the uniform distribution of the 4-coupled symmetric tent map is plotted on Figure 9 and its projection on the phase subspace $(x^1, x^3)$ is displayed on Figure 10.

$N_{\text{disc}}$ is fixed to $10^2 \times 10^2$ and $N_{\text{iter}}$ to $10^{11}$, $\varepsilon_1$ to $10^{-14}$.

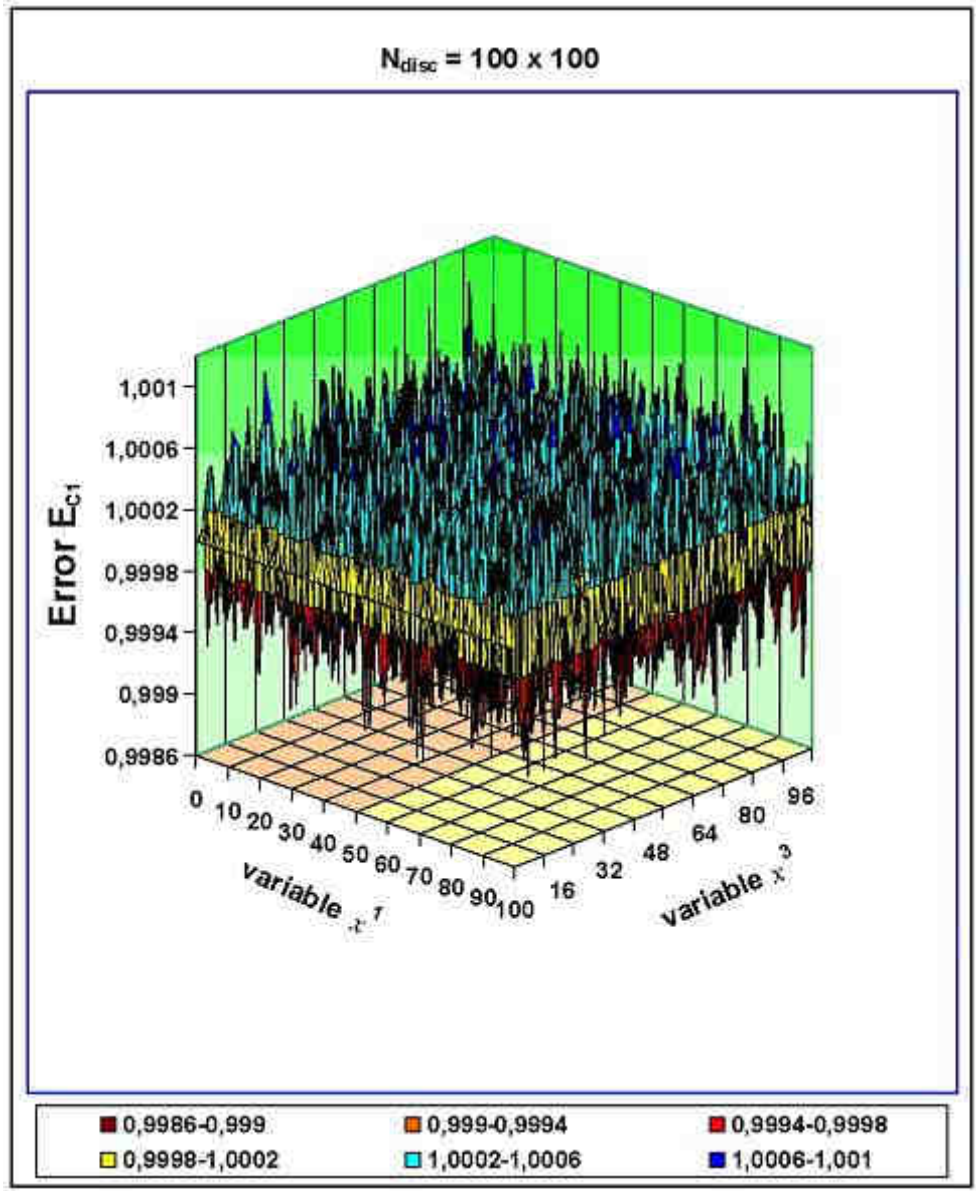

**Figure 9.** Difference between the correlation distribution function $C_{N_{DISC}}(x^1, x^3)$ and the uniform distribution of the 4-coupled symmetric tent map. $N_{disc} = 10^2 \times 10^2$, $N_{iter} = 10^{11}$, $\varepsilon_i = i.\varepsilon_1$, $\varepsilon_1 = 10^{-14}$. Computations done using double precision numbers(~14-15 digits). Initial values $x^1_0 = 0.330$, $x^2_0 = 0.3387564$, $x^3_0 = 0.50492331$, $x^4_0 = 0.0$.

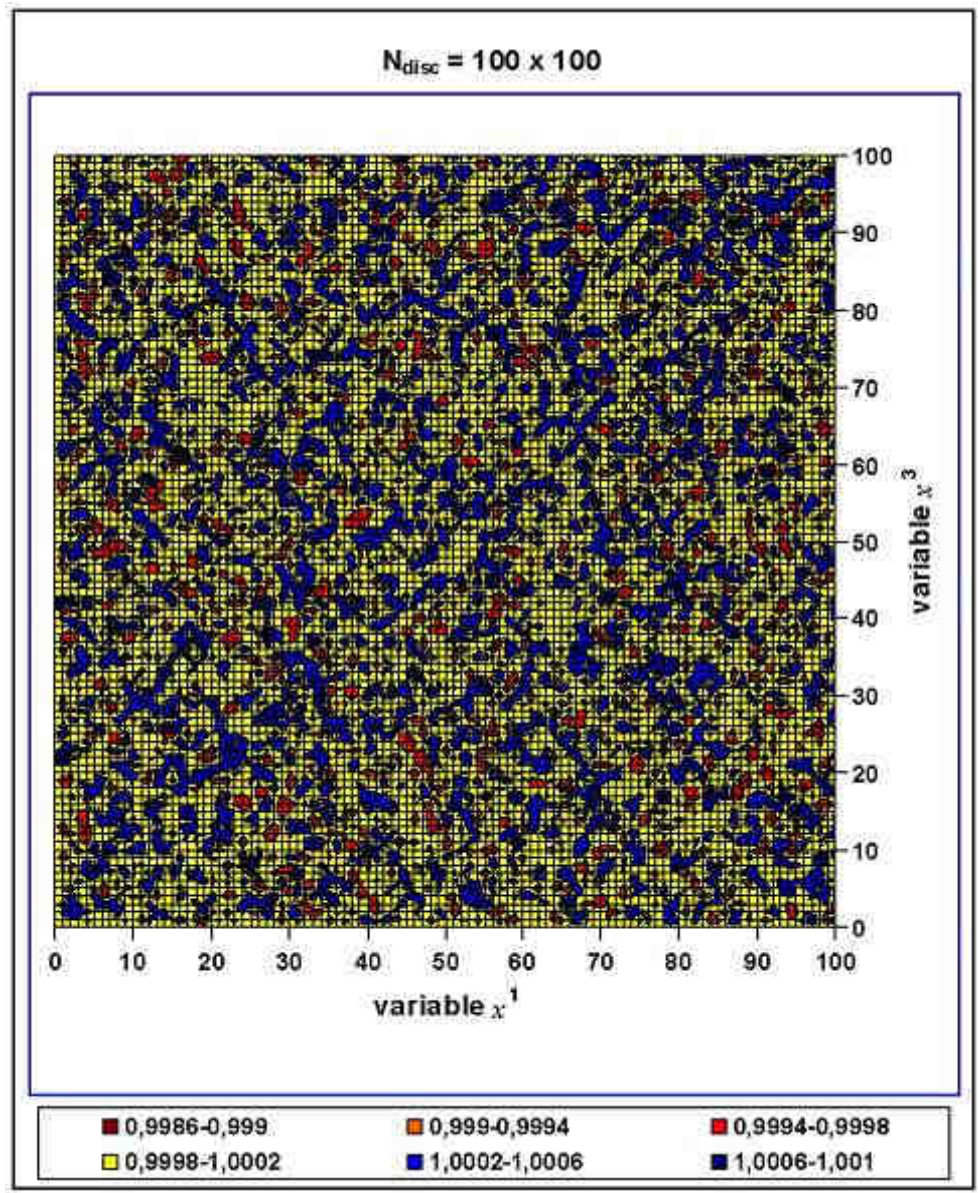

**Figure 10.** Projection of the difference between the correlation distribution function $C_{N_{DISC}}(x^1, x^3)$ and the uniform distribution of the 4-coupled symmetric tent map on the phase subspace $(x^1, x^3)$. $N_{disc} = 10^2 \times 10^2$, $N_{iter} = 10^{11}$, $\varepsilon_i = i.\varepsilon_1$, $\varepsilon_1 = 10^{-14}$. Computations done using double precision numbers(~14-15 digits). Initial values $x^1_0 = 0.330$, $x^2_0 = 0.3387564$, $x^3_0 = 0.50492331$, $x^4_0 = 0.0$.

In order to check the uncorrelation, every couple of components must be checked simultaneously. In the considered case $N_{iter} = 10^{11}$ for the 4-coupled symmetric tent map, the errors $E_{c1(N_{disc}, N_{iter})}(x^k, x^l)$ for $1 \leq k \leq l \leq 4$ are displayed in Table 7.

| $E_{c1(N_{disc},N_{iter})}(x^k,x^l)$ | $x^2$ | $x^3$ | $x^4$ |
|---|---|---|---|
| $x^1$ | 0.00025640111 | 0.00025509367 | 0.00025269625 |
| $x^2$ |  | 0.00025215219 | 0.00025067251 |
| $x^3$ |  |  | 0.00024861781 |

**Table 7.** Error $E_{c1(N_{disc},N_{iter})}(x^k,x^l)$ for $1 \leq k \leq l \leq 4$.

$N_{disc} = 10^2 \times 10^2$, $N_{iter} = 10^{11}$, $\varepsilon_i = i.\varepsilon_1$, $\varepsilon_1 = 10^{-14}$. Computations done using double precision numbers(~14-15 digits). Initial values $x^1_0 = 0.330$, $x^2_0 = 0.3387564$, $x^3_0 = 0.50492331$, $x^4_0 = 0.0$.

## 3. Enhanced chaotic generators

### 3.1 Chaotic sampling of chaotic numbers

If we plot the chaotic numbers produced by any component $x^l$, $1 \leq l \leq p$ of the $p$-dimensional dynamical system (2.12) in the phase space $(x_n^l, x_{n+1}^l)$, the iterated points show up the graph of the symmetrical tent map $f$ used to define (2.12) (more exactly a graph with two lines of $\varepsilon$ thickness). These numbers are not randomly produced. If we plot these points in the phase space $(x_n^l, x_{n+2}^l)$ or $(x_n^l, x_{n+3}^l)$ they will display the graph of $f^{(2)}$ or $f^{(3)}$ (see Figure 11). Hence, someone knowing a sequence of few iterated points is able to find the initial value $X_0$ of the dynamical system.

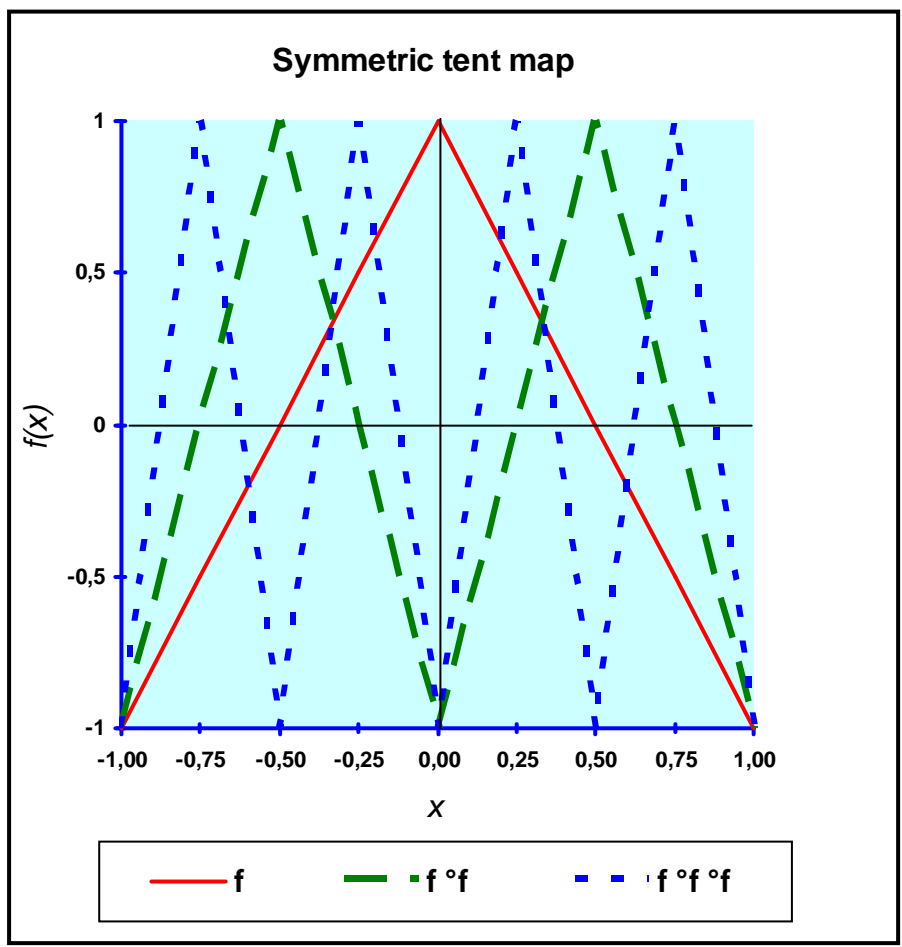

**Figure 11.** Graphs of the symmetric tent map $f$, $f^{(2)}$ and $f^{(3)}$ on the interval [-1,1].

In order to hide the graph of the genuine function $f$ in the phase space $(x_n^l, x_{n+1}^l)$, a pivotal idea is to sample chaotically the sequence $(x_0^l, x_1^l, x_2^l, \ldots, x_n^l, x_{n+1}^l, \ldots)$ selecting $x_n^l$ every time the value of $x_n^m$ is greater than a threshold $T$, $-1 < T < 1$, with $l \neq m$.

The chaotically sampled subsequence $(\overline{x_0}, \overline{x_1}, \overline{x_2}, \cdots, \overline{x_q}, \overline{x_{q+1}}, \cdots)$ being defined as:

$$\overline{x_q} = x_n^l \quad iff \quad x_n^m \in\ ]T, 1[ \quad (3.2)$$

Choosing $T > 0.5$ implies that the selected subsequence

$$(\overline{x_0}, \overline{x_1}, \overline{x_2}, \cdots, \overline{x_q}, \overline{x_{q+1}}, \cdots) = (x_{p_0}^l, x_{p_1}^l, x_{p_2}^l, \ldots, x_{p_q}^l, x_{p_{q+1}}^l, \cdots)$$

is such that the difference between $p_q$ and $p_{q+1}$ is always greater than a minimal value $K_m$ depending upon $T$.

The graph of the chaotically sampled chaotic number is a mixing of the graphs of all the $f^{(r)}$ for $r > K_m$.

Since we have shown in section 2.2.4 that every pair of components $(x_n^l, x_n^m)$ of the sequence $X_0, X_1, X_2, \cdots, X_n, X_{n+1}, \cdots$ is uncorrelated, the chaotic sampling process we propose, is a powerful tool to generate enhanced chaotic numbers. We perform some numerical experiments in this section in order to check their properties.

We define $AC_{M,N}(x, y)$ the autocorrelation distribution function which is the correlation function $C_{M,N}(x, y)$ (see 2.2.4) in the phase space $(x_n^l, x_{n+1}^l)$ instead of the phase space $(x^l, x^m)$. In order to control that the enhanced chaotic numbers $(\overline{x_0}, \overline{x_1}, \overline{x_2}, \cdots, \overline{x_q}, \overline{x_{q+1}}, \cdots)$ are uncorrelated, we plot them in the phase subspace $(\overline{x_n}, \overline{x_{n+1}})$ and we check if they are uniformly distributed in the square $[-1, 1]^2$.

If no particular pattern is displayed and if the difference between the autocorrelation distribution function $AC_N(x, y)$ converges towards the uniform distribution on the square when the number of iterations goes to the infinity, we can conclude that the knowledge of a sequence of iterated points do not allow finding the initial value $X_0$ of the dynamical system.

On Figure 12 are displayed the values of

$$E_{AC1}(N_{disc}, NSampl_{iter}) = \left\| AC_{N_{disc}, NSampl_{iter}}(x, y) - 0.25 \right\|_{L_1}$$ for a system of 4 coupled-equations for both threshold values 0.98 and 0.998 of $x_n^4$. The enhanced chaotic numbers are produced by the first component $x_n^1$ of the dynamical system. As the chaotic numbers are regularly distributed on

the interval [-1, 1], when $T > 0.98$ one chaotic number over approximately 100 is sampled, when $T > 0.998$ one chaotic number over approximately 1,000 is sampled. We call $NSampl_{iter}$ the number of sampled points.

Nevertheless the computing process is very fast. A desktop computer can produce more than 50,000,000 chaotic numbers per second, thus 50,000 iterated sampled points per second for $T > 0.998$. The sampling threshold 0.998 gives very good results.

The difference between the autocorrelation distribution function $AC_{\text{NSAMPLDISC}}\left(\overline{x_n}, \overline{x_{n+1}}\right)$ and the uniform distribution is plotted on Figure 13 and its projection on the phase subspace $\left(\overline{x_n}, \overline{x_{n+1}}\right)$ is displayed on Figure 14.

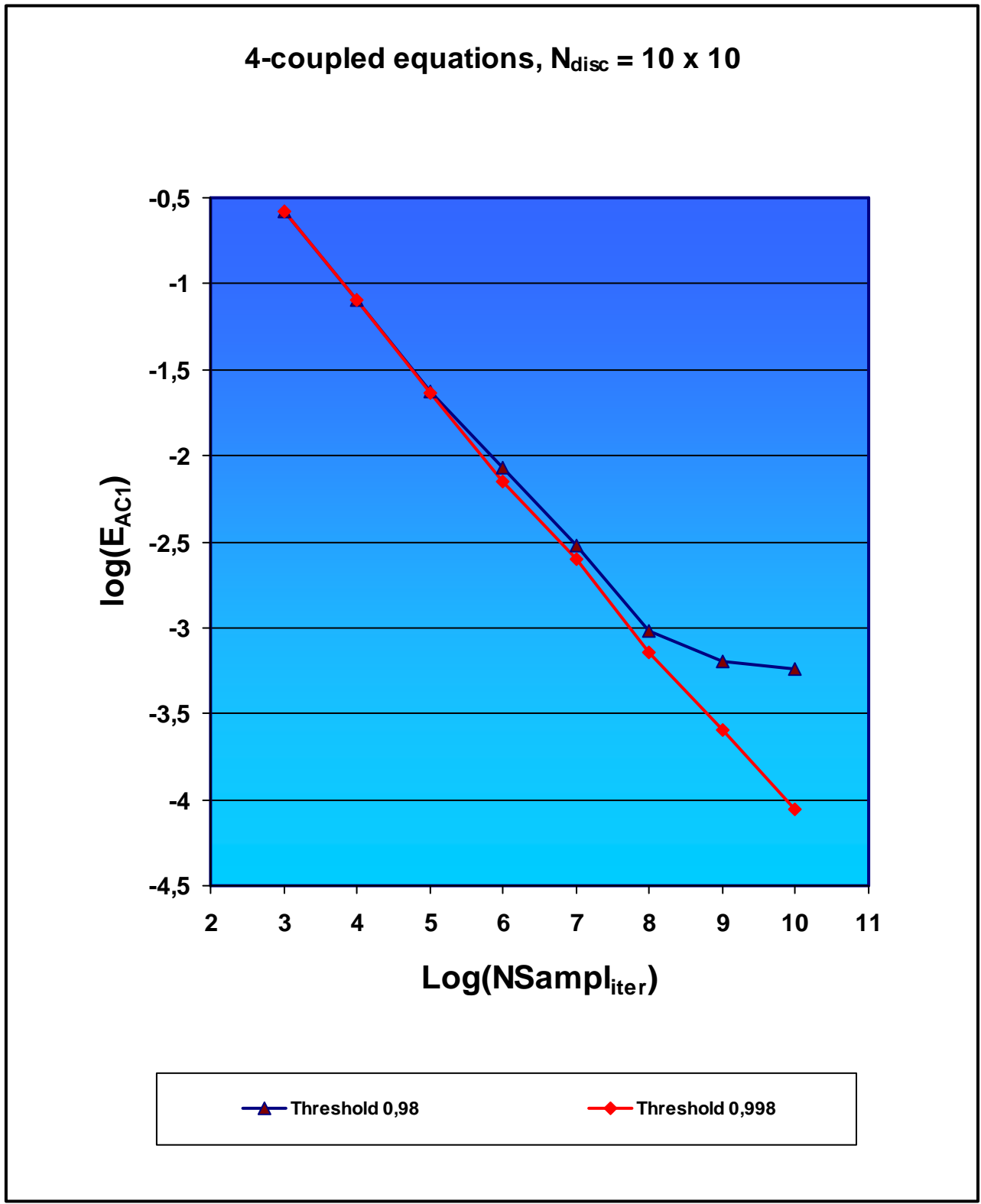

**Figure 12.** Error $E_{AC1}$ for the first component $x^1$, sampled by $x^4$ for the threshold values 0.98 and 0.998 of the 4-coupled symmetric tent map. $N_{disc} = 10 \times 10$, $\varepsilon_i = i.\varepsilon_1$, $\varepsilon_1 = 10^{-14}$, $NSampl_{iter}$ varies from $10^3$ to $10^{10}$. Computations done using double precision numbers (~14-15 digits). Initial values: $x^1_0 = 0.330$, $x^2_0 = 0.3387564$, $x^3_0 = 0.50492331$, $x^4_0 = 0.0$.

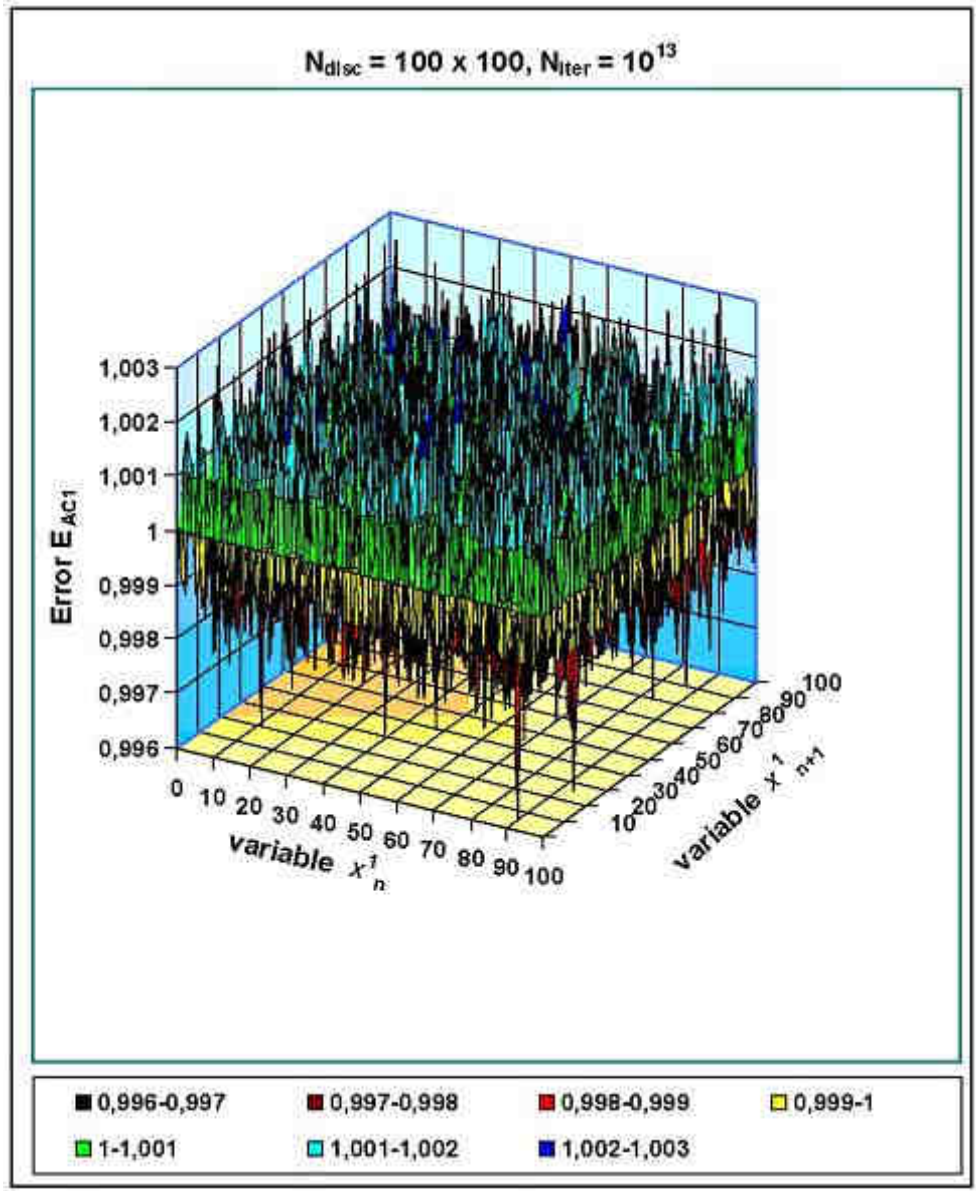

**Figure 13.** Difference between the autocorrelation distribution function $AC_{\text{NSAMPLDISC}}(x_n^1, x_{n+1}^1)$ and the uniform distribution of the 4-coupled symmetric tent map sampled by $x^4$ for the threshold values 0.98 and 0.998 of. $N_{\text{disc}} = 10^2 \times 10^2$, $NSampl_{\text{iter}} = 10^{10}$, $\varepsilon_i = i.\varepsilon_1$, $\varepsilon_1 = 10^{-14}$. Initial values $x^1_0 = 0.330$, $x^2_0 = 0.3387564$, $x^3_0 = 0.50492331$, $x^4_0 = 0.0$.

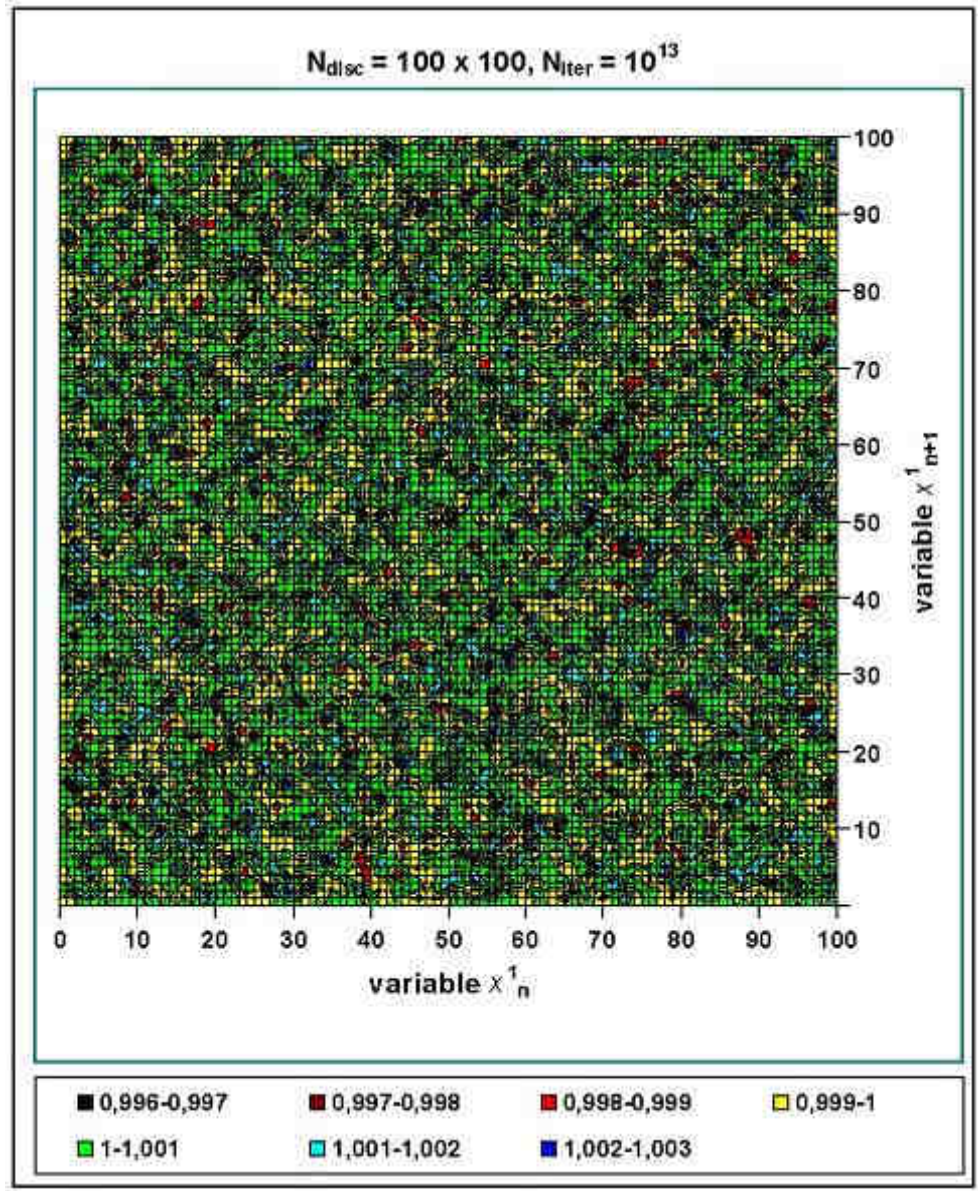

**Figure 14.** Projection of the difference between the autocorrelation distribution function $AC_{\text{NSAMPLDISC}}\left(x_n^1, x_{n+1}^1\right)$ and the uniform distribution of the 4-coupled symmetric tent map sampled by $x^4$ on the phase subspace $\left(x_n^1, x_{n+1}^1\right)$. $N_{\text{disc}} = 10^2 \times 10^2$, $NSampl_{\text{iter}} = 10^{10}$, $\varepsilon_i = i.\varepsilon_1$, $\varepsilon_1 = 10^{-14}$. Initial values $x^1_0 = 0.330$, $x^2_0 = 0.3387564$, $x^3_0 = 0.50492331$, $x^4_0 = 0.0$.

Figure 15 and Table 8 display $E_{AC1}(N_{disc}, NSampl_{iter})$ with respect to $N_{\text{disc}}$.

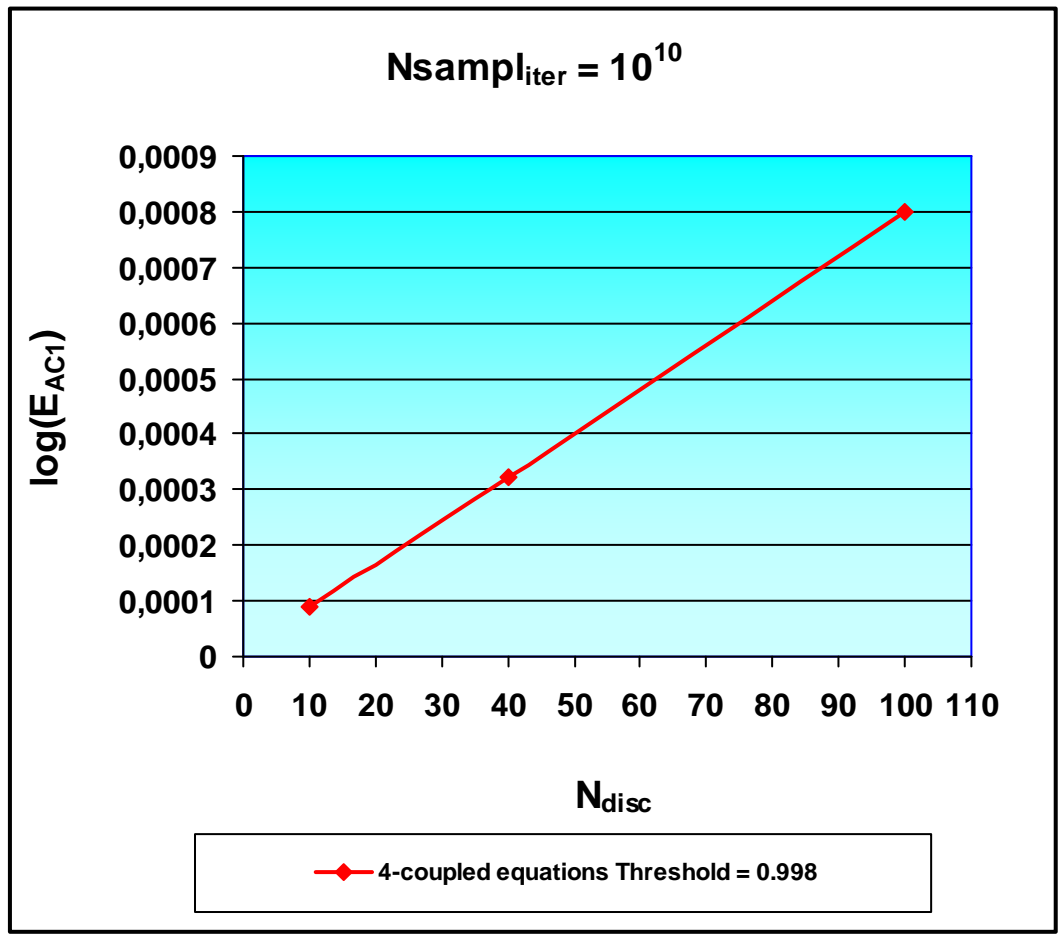

**Figure 15.** $E_{AC1}(N_{disc}, NSampl_{iter})$ for the first component $x^1$, sampled by $x^4$ for the threshold value 0.998 of the 4-coupled symmetric tent map versus $N_{disc}$, $NSampl_{iter} = 10^{10}$, $\varepsilon_i = i.\varepsilon_1$, $\varepsilon_1 = 10^{-14}$. Initial values $x^1_0 = 0.330$, $x^2_0 = 0.3387564$, $x^3_0 = 0.50492331$, $x^4_0 = 0.0$.

| $N_{disc}$ | $NSampl_{iter}$ | $E_{AC1}(N_{disc}, NSampl_{iter})$ |
|---|---|---|
| **10 x 10** | **10,000,042,552** | 0.000088445108 |
| **40 x 40** | **10,000,042,552** | 0.00032254866 |
| **100 x 100** | **10,000,042,552** | 0.00079801406 |

**Table 8.** $E_{AC1}(N_{disc}, NSampl_{iter})$ for the first component $x^1$, sampled by $x^4$ for the threshold value 0.998 of the 4-coupled symmetric tent map versus $N_{disc}$, $NSampl_{iter} = 10^{10}$, $\varepsilon_i = i.\varepsilon_1$, $\varepsilon_1 = 10^{-14}$. Initial values $x^1_0 = 0.330$, $x^2_0 = 0.3387564$, $x^3_0 = 0.50492331$, $x^4_0 = 0.0$.

### 3.1 Mixing and chaotic sampling of chaotic numbers

One can improve the unpredictability of the chaotic numbers generated as above, using all the components of the vector $X = \begin{pmatrix} x^1 \\ \vdots \\ \vdots \\ x^p \end{pmatrix}$ instead of one. For example for 4-coupled equations, the value of $x_n^4$ command the sampling process in the following way:

Let set three threshold values $T_1$, $T_2$ and $T_3$

$$-1 < T_1 < T_2 < T_3 < 1 \qquad (3.2)$$

we sample and mix together chaotically the sequences $(x_0^1, x_1^1, x_2^1, \ldots, x_n^1, x_{n+1}^1, \ldots)$, $(x_0^2, x_1^2, x_2^2, \ldots, x_n^2, x_{n+1}^2, \ldots)$ and $(x_0^3, x_1^3, x_2^3, \ldots, x_n^3, x_{n+1}^3, \ldots)$ defining $(\overline{x_0}, \overline{x_1}, \overline{x_2}, \cdots, \overline{x_q}, \overline{x_{q+1}}, \ldots)$ by:

$$\overline{x_q} = \begin{cases} x_n^1 & \text{iff} \quad x_n^4 \in \, ]T_1, T_2[ \\ x_n^2 & \text{iff} \quad x_n^4 \in \, [T_2, T[ \\ x_n^3 & \text{iff} \quad x_n^4 \in \, [T_3, 1[ \end{cases} \qquad (3.3)$$

On Figure 16 and Table 9 are displayed the values of $E_{AC1}(N_{disc}, NSampl_{iter})$ for a system of 4 coupled-equations when the first component $x^1$ is sampled by $x^4$ for both the threshold values 0.98 and 0.998 and when the three components $x^1$, $x^2$, $x^3$ are mixed and sampled by $x^4$ for the threshold values $T_1 = 0.98$, $T_2 = 0.987$, $T_3 = 0.994$ or $T_1 = 0.998$, $T_2 = 0.9987$, $T_3 = 0.9994$.

On Figure 17, in order to visualize the converge of the autocorrelation function towards the uniform repartition we have plotted the difference between the autocorrelation distribution function $AC_{\text{NSAMPLDISC}}(\overline{x_n}, \overline{x_{n+1}})$ and the uniform distribution of the 4-coupled symmetric tent when the three components $x^1$, $x^2$, $x^3$ are mixed and sampled by $x^4$ for the threshold values $T_1 = 0.998$, $T_2 = 0.9987$, $T_3 = 0.9994$. $NSampl_{\text{iter}} = 10^9$ on the left part of the figure and $10^{10}$ on the right.

## 4. Conclusion

We have introduced new families of enhanced chaotic number generators in order to hide the generating function. The key feature of these enhanced chaotic number generators is that they use chaotic numbers themselves in order to sample and mix chaotic subsequences of chaotic numbers. The properties of these new families are explored numerically up to $10^{13}$ iterations.

The numerical experiments give good results; other tests have to be performed in order to check their usefulness as Chaotic PRNG.

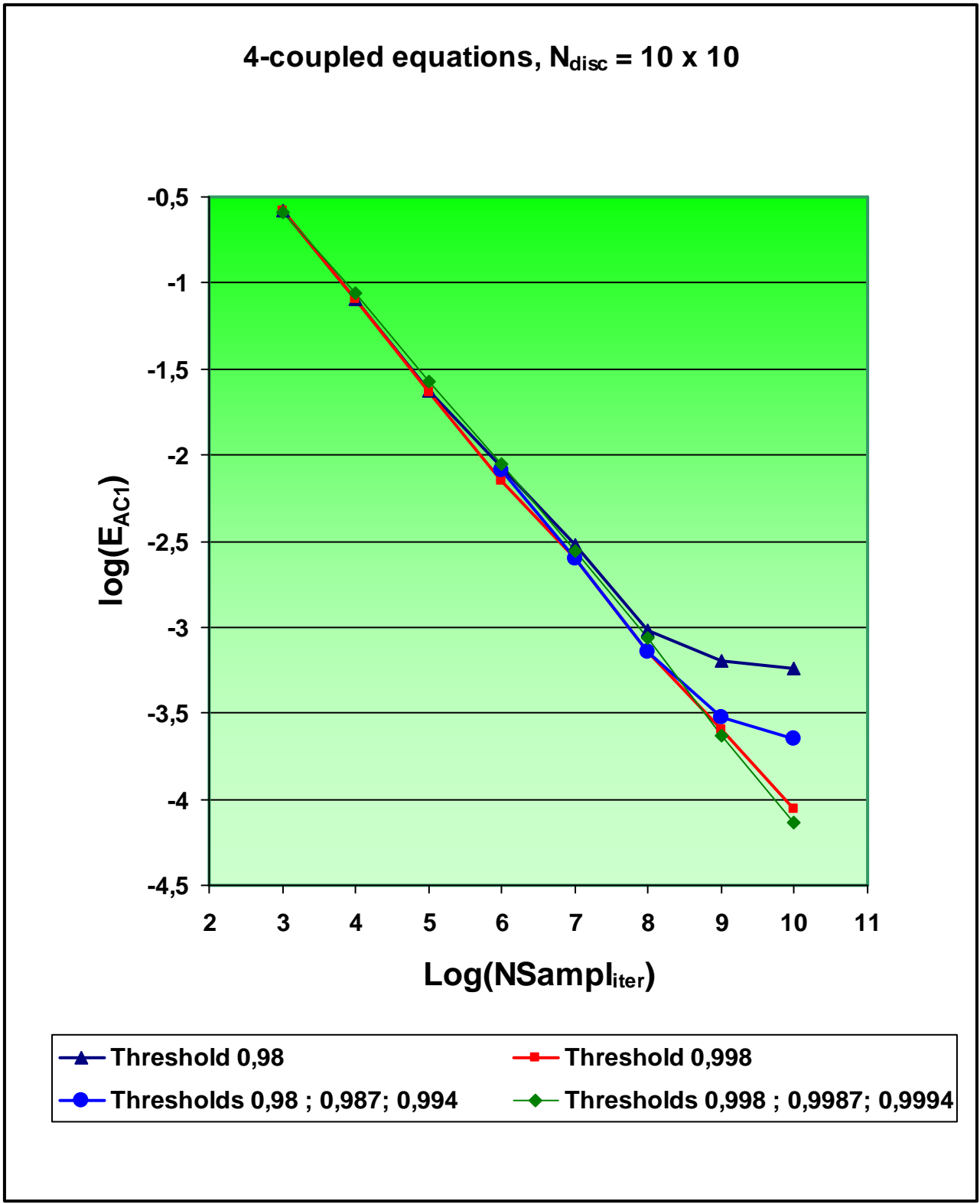

**Figure 16.** Error of $E_{AC1}(N_{disc}, NSampl_{iter})$ for a system of 4 coupled-equations when the first component $x^1$ is sampled by $x^4$ for both the threshold values 0.98 and 0.998 and when the three components $x^1$, $x^2$, $x^3$ are mixed and sampled by $x^4$ for the threshold values $T_1 = 0.98$, $T_2 = 0.987$, $T_3 = 0.994$ and $T_1 = 0.998$, $T_2 = 0.9987$, $T_3 = 0.9994$. $N_{disc} = 10 \times 10$, $\varepsilon_i = i.\varepsilon_1$, $\varepsilon_1 = 10^{-14}$, $NSampl_{iter}$ varies from $10^3$ to $10^{10}$. Initial values: $x^1_0 = 0.330$, $x^2_0 = 0.3387564$, $x^3_0 = 0.50492331$, $x^4_0 = 0.0$.

| $N_{\text{iter}}$ | $NSampl_{\text{iter}}$ | $E_{AC1}(N_{disc}, NSampl_{iter})$<br>**4-coupled equation**<br>$T = 0.998$ | $NSampl_{\text{iter}}$ | $E_{AC1}(N_{disc}, NSampl_{iter})$<br>**4-coupled equation**<br>$T_1 = 0.998, T_2 = 0.9987,$<br>$T_3 = 0.9994$ |
|---|---|---|---|---|
| $10^5$ | 95 | 0.70947368 | 93 | 0.68924731 |
| $10^6$ | 971 | 0.26570546 | 1015 | 0.25881773 |
| $10^7$ | 10,095 | 0.079871223 | 10,139 | 0.086706776 |
| $10^8$ | 100,622 | 0.023190157 | 100,465 | 0.026815309 |
| $10^9$ | 1,001,408 | 0.0071386288 | 1,000,549 | 0.0089111078 |
| $10^{10}$ | 9,998,496 | 0.002493667 | 9,998,814 | 0.0027932033 |
| $10^{11}$ | 100,013,867 | 0.00071561417 | 100,001,892 | 0.00085967214 |
| $10^{12}$ | 999,994,003 | 0.00025442753 | 999,945,728 | 0.0002346851 |
| $10^{13}$ | 10,000,042,552 | 0.000088445108 | 10,000,046,137 | 0.000073234736 |

**Table 9.** Error of $E_{AC1}(N_{disc}, NSampl_{iter})$ for a system of 4 coupled-equations when the first component $x^1$ is sampled by $x^4$ for both the threshold values 0.98 and 0.998 and when the three components $x^1$, $x^2$, $x^3$ are mixed and sampled by $x^4$ for the threshold values $T_1 = 0.98$, $T_2 = 0.987$, $T_3 = 0.994$ and $T_1 = 0.998$, $T_2 = 0.9987$, $T_3 = 0.9994$. $N_{disc} = 10 \times 10$, $\varepsilon_i = i.\varepsilon_1$, $\varepsilon_1 = 10^{-14}$, $NSampl_{\text{iter}}$ varies from $10^3$ to $10^{10}$. Initial values: $x^1_0 = 0.330$, $x^2_0 = 0.3387564$, $x^3_0 = 0.50492331$, $x^4_0 = 0.0$.

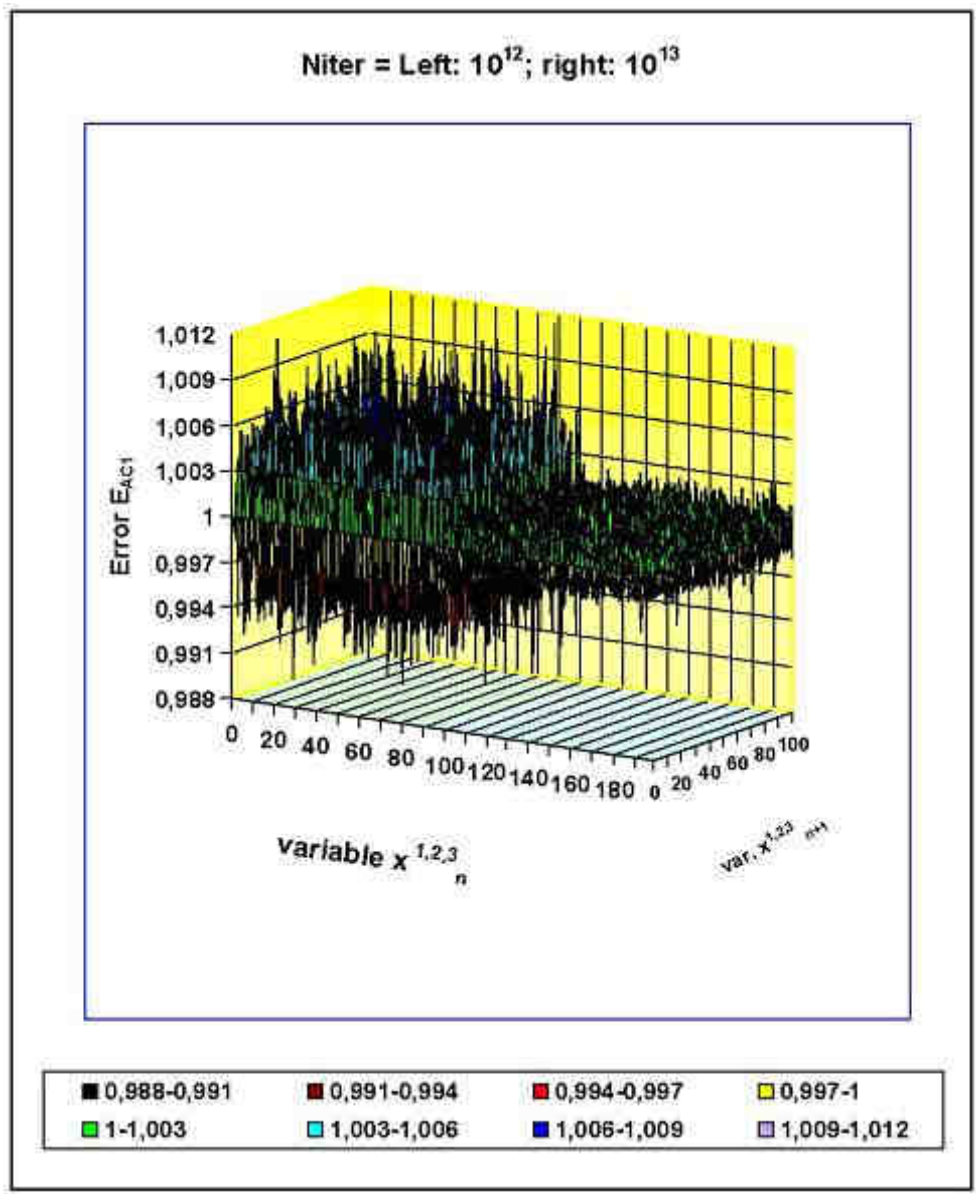

**Figure 17.** Difference between the autocorrelation distribution function $AC_{\text{NSAMPLDISC}}(x_n^1, x_{n+1}^1)$ and the uniform distribution of the 4-coupled symmetric tent when the three components $x^1$, $x^2$, $x^3$ are mixed and sampled by $x^4$ for the threshold values $T_1 = 0.998$, $T_2 = 0.9987$, $T_3 = 0.9994$. $N_{\text{disc}} = 100 \times 100$, $\varepsilon_i = i.\varepsilon_1$, $\varepsilon_1 = 10^{-14}$, $NSampl_{\text{iter}} = 10^9$ on the left part of the figure and $10^{10}$ on the right. Initial values: $x^1_0 = 0.330$, $x^2_0 = 0.3387564$, $x^3_0 = 0.50492331$, $x^4_0 = 0.0$.